\documentclass[12pt]{article}
\usepackage{amssymb,amsmath,graphicx}
\usepackage{theorem}
\usepackage[usenames]{color}
\usepackage{hyperref}
\date{May 2006\\(revised January 2007)}
\nonstopmode

\newtheorem{theorem}{Theorem}[section]
\newtheorem{prop}[theorem]{Proposition}
\newtheorem{lemma}[theorem]{Lemma}
\newtheorem{claim}[theorem]{Claim}
\newtheorem{corollary}[theorem]{Corollary}
\newtheorem{suppl}[theorem]{Supplement}

\theorembodyfont{\rmfamily}

\newtheorem{definition}[theorem]{Definition}
\newtheorem{example}{Example}[section]

\newtheorem{remark}[theorem]{Remark}

\numberwithin{equation}{section}

\newenvironment{proof}{\medskip\noindent{\bf Proof.
}}{\hfill$\square$\medskip}

\def\bsa{\begin{subarray}{c} }
\def\esa{\end{subarray} }

\def\hom{{\rm hom}}
\def\inj{{\rm inj}}

\def\Inj{{\rm Inj}}

\def\inj{{\rm inj}}
\def\ind{{\rm ind}}

\def\eps{\varepsilon}

\def\IO{{\infty\to1}}

\def\iv{\widehat}

\def\Prob{{\sf P}}
\def\E{{\sf E}}
\def\Var{{\sf Var}}

\def\T{{^\top}}

\def\one{{\mathbf 1}}

\def\bG{{\mathbf G}}
\def\bH{{\mathbf H}}
\def\bX{{\mathbf X}}

\def\RandG{{\mathbb G}}

\def\ba{{\mathbf a}}

\def\R{{\mathbb R}}

\def\Z{{\mathbb Z}}

\def\AA{{\cal A}}

\def\MM{{\cal M}}
\def\PP{{\cal P}}\def\QQ{{\cal Q}}

\def\WW{{\cal W}}\def\XX{{\cal X}}

\def\WWz{{\WW_{[0,1]}}}

\def\ontop#1#2{\genfrac{}{}{0pt}{}{#1}{#2}}

\begin{document}

\title{Convergent Sequences of Dense Graphs I: \\
Subgraph Frequencies,
Metric Properties and Testing}

\author{C.~Borgs$^{a}$,\;
J.T.~Chayes$^a$,\;
L.~Lov\'asz$^a$%
\thanks{Current Address: E\"otv\"os Lor\'and University,
P\'azm\'any P\'eter s\'et\'any 1/C, H-1117 Budapest,
Hungary},\;
V.T.~S\'os$^b$%
\thanks{Research supported in part by OTKA grants T032236,
T038210, T042750},\;
\and
K.~Vesztergombi$^c$
\\
\\
$^a${\small\it Microsoft Research, One Microsoft Way, Redmond,
WA 98052, USA}
\vspace{-6pt}\and
$^b${\small\it Alfr\'ed R\'enyi Institute of Mathematics,
POB 127, H-1364 Budapest, Hungary}
\vspace{-5pt}\and
$^c${\small\it E\"otv\"os Lor\'and University,
P\'azm\'any P\'eter S\'etani 1/C, H-1117 Budapest,
Hungary}\\[5pt]}
\maketitle

\begin{abstract}

We consider sequences of graphs $(G_n)$ and define various notions of
convergence related to these sequences: ``left convergence'' defined
in terms of the densities of homomorphisms from small graphs into
$G_n$; ``right convergence'' defined in terms of the densities of
homomorphisms from $G_n$ into small graphs; and convergence in
a suitably defined metric.

In Part I of this series, we show that left convergence is equivalent
to convergence in metric, both for simple graphs $G_n$, and for
graphs $G_n$ with nodeweights and edgeweights.
One of the main steps here is the introduction of a cut-distance
comparing graphs, not necessarily of the same size.
We also
show how these notions of convergence provide natural formulations of
Szemer\'edi partitions, sampling and testing of large graphs.
\end{abstract}


\tableofcontents

\addtolength{\baselineskip}{3pt}
\section{Introduction}
\label{sec:intro}

In this and accompanying papers, we define a natural notion
of convergence of a sequence of graphs, and show that other useful
notions of convergence are equivalent to it.  We are motivated by
the fact that, in many subfields of
mathematics, computer science and physics, one studies
properties of very large graphs, or properties of graph sequences
that grow beyond all limits.  Let us give a few examples:

\smallskip

{\it Random Networks:} There is a large literature of graph models
of the Internet, the WWW, and other so-called ''scale-free''
technological networks, {first modeled in this context by
Barab\'asi and Albert \cite{BA}. Among the technological networks
that are modeled are} the graph of computers and physical links
between them, the graph of {so-called Autonomous Systems such as
Internet Service Providers,} the graph of webpages with hyperlinks,
etc. These graphs are often similar to various graphs of social
networks: acquaintances, co-publications, the spreading of certain
diseases, etc.

These networks are formed by random processes, but their properties
are quite different from the traditional Erd\H{o}s--R\'enyi random
graphs: their degree distribution has a ``heavy tail'', they tend to
be clustered, the neighborhoods of their nodes are denser than the
average edge density, etc. Several models of random {scale-free}
graphs have been proposed and studied. {For rigorous work, see
\cite{BR,BRST} for undirected models, \cite{copy} for ``copying
models'', \cite{BBCR} for a directed model, \cite{BBCS} for the
spread of viruses on these networks, and \cite{BR-rev} for a survey
of rigorous work with more complete references.}

{\it Quasirandom Graphs:  } Quasirandom (also called pseudorandom)
graphs were introduced by Thomason \cite{Tho} and Chung, Graham and
Wilson \cite{CGW}. These graph sequences {can be} deterministic, but
have many properties of true random graphs. A nice example is the
sequence of Paley graphs (quadratic residue graphs).
These graphs are remarkably similar to a random graph with
edge-probability $1/2$ on the same number of nodes in many ways. The
most relevant for us is that they contain (asymptotically) the same
number of copies of each fixed graph $F$ as the random
graph---this is one of the many equivalent ways to define
quasirandom graphs. Many other questions in graph theory, in
particular in extremal graph theory, also involve asymptotic counting
of small graphs.

{\it Property Testing of Large Graphs:} Say we are given a large
graph and we want to determine certain numerical parameters, e.g.,
the edge density, of that graph by sampling a bounded number of
nodes. Or perhaps we want to determine whether the large graph has a
given property, e.g., is it 3-colorable?  In particular, which
parameters can be accurately estimated and which properties can be
tested with high probability by looking only at subgraphs on small
randomly chosen subsets of the nodes? A precise definition of
property testing was given by Goldreich, Goldwasser and Ron
\cite{GGR}, who also proved several fundamental results about this
problem.

{\it Statistical Mechanics:} Many models in physics are described by
a weighted coloring of some large graph $G$.  The graph $G$ typically
represents underlying geometric structure of the model under
consideration, e.g. a crystal lattice and its nearest neighbor
structure, while the color of a given node represents the local
state. In the simplest case of two colors, the two colors can
represent quantities like the two possible orientations of a spin
variable, or the presence or absence of a molecule at a given
position.  The interactions between different local states can then
be described by a weighted ``interaction graph'' $H$, with smaller
edgeweights corresponding to weaker interactions, and larger
edgeweights representing stronger interactions. In this context, the
weighted number of colorings represents the so-called partition
function of the model.

{\it Combinatorial Optimization:} Many optimization problems can be
described as weighted coloring problems. A simple example is the
max-cut problem, where our task is to find the maximal cut in a large
graph $G$. If we consider a coloring of $G$ with two colors, $1$ and
$2$, and weight a coloring by the number of edges with two
differently colored endnodes, then
the maximum cut is just given by the maximum weight coloring.

\medskip

In this and two accompanying papers we develop a theory of
convergence of  graph sequences, which works best in two extreme
cases: dense graphs (the subject of this paper and \cite{dense2}) and
graphs with bounded degree (the subject of \cite{sparse}).
Convergence of graph sequences was defined by Benjamini and Schramm
\cite{BS} for graphs with bounded degree, and by the authors of this
paper \cite{BCLSV-03} for dense graphs.

Our general setup will be the following. We have a ``large'' graph
$G$ with node set $V(G)$ and edge set $E(G)$. There are (at least)
two ways of studying $G$ using homomorphisms. First, we can count the
number of copies of various ``small'' graphs $F$ in $G$, more
precisely, we count the number of homomorphisms from $F$ to $G$; this
way of looking at $G$ allows us to treat many problems in, {e.g.,}
extremal graph theory. Second, we can count homomorphisms from
$G$ into various small graphs $H$; this includes
many models in
statistical physics and many problems on graph coloring.

These two notions of probing a large graph with a small graph lead
to two different notions of convergence of a sequence of graphs
$(G_n)$: convergence from the left, corresponding to graphs which
look more and more similar when probed with homomorphisms from small
graphs into $G_n$, and convergence from the right, corresponding to
graph sequences whose elements look more and more similar when
probed with homomorphism from $G_n$ into a small graphs.

This theory can also be viewed as a substantial generalization of the
theory of quasirandom graphs. In fact, most of the equivalent
characterizations of quasirandom graphs are immediate corollaries of
the general theory developed here and in our companion paper
\cite{dense2}.

In this paper we study convergence from the left, both for sequences
of simple graphs and sequences of weighted graphs, and its relations
to sampling and testing.  {Since this paper focuses on convergence
from the left, we'll often omit the phrase "from the left".} We will
also show that convergence from the left is equivalent to convergence
in metric for a suitable
notion of distance between two weighted graphs. Finally, we
will show that convergence from the left is equivalent to the
property that the graphs in the sequence have asymptotically the same
Szemer\'edi partitions.

Convergence from the right will be the subject matter of the sequel
of this paper \cite{dense2}.

Convergence in metric clearly allows for a completion by the usual
abstract identification of Cauchy sequences of distance zero.  But it
turns out (Lov\'asz and Szegedy \cite{LSz1}) that the limit
object of a convergent graph sequence has a much more natural
representation in terms of a measurable symmetric function
$W:[0,1]^2\to\R$ (we call these functions {\it graphons}). In
fact, it is often useful to represent a finite graph $G$ in terms of
a suitable function $W_G$ on $[0,1]^2$, defined as step
function with steps of length $1/|V(G)|$ and
values $0$ and $1$, see below for the precise definition.
While the introduction of graphons requires some basic notions of
measure theory, it will simplify many proofs in this paper.

The organization of this paper is as follows: In the next section, we
introduce our definitions: in addition to left-convergence, we define
a suitable distance between weighted graphs, and state our main
results for weighted graphs.  In Section~\ref{sec:graphon}, we
generalize these definitions and results to graphons. The following
section, Section~\ref{sec:sampling}, is devoted to sampling, and
contains the proofs of the main results of this paper, including the
equivalence of left-convergence and convergence in metric.
Section~\ref{sec:norm+szem} relates convergence in metric to an
{\it a
priori} weaker form of ``convergence in norm'' and
to Szemer\'edi partitions, and Section~\ref{sec:testing} proves
our results on testing.  We close this paper with a section on
miscellaneous results and an outlook on right-convergence. In
the appendix, we describe a few details of proofs
which are omitted in the main body of the paper.

\section{Weighted and Unweighted Graphs}
\label{sec:graphs}

\subsection{Notation}
\label{sec:Notation}

We consider both unweighted, simple graphs and weighted graphs,
where, as usual, a simple graph $G$ is a graph without loops or
multiple edges. We denote the node and edge set of $G$ by $V(G)$ and
$E(G)$, respectively.

A {\it weighted graph} $G$ is a graph with a weight
$\alpha_i=\alpha_i(G)>0$ associated with each node and a weight
$\beta_{ij}=\beta_{ij}{(G)}\in\R$ associated with each edge $ij$,
including possible loops with $i=j$. For convenience, we set
$\beta_{ij}=0$ if $ij\notin E(G)$. We set
\[\alpha_G=\sum_i\alpha_i(G),\quad
\|G\|_\infty=\max_{i,j}|\beta_{ij}(G)|,
\quad\text{and}\quad
\|G\|_2=
\Biggl(\sum_{i,j}\frac{\alpha_i\alpha_j}{\alpha_G^2}\beta_{ij}^2
\biggr)^{1/2},
\]
and for $S,T\subset V(G)$, we define
\begin{equation}
\label{eG-def}
e_G(S,T)=\sum_{\bsa i\in S\\ j\in T\esa}
\alpha_i(G)\alpha_j(G)\beta_{ij}(G).
\end{equation}

A weighted graph $G$ is called {\it soft-core} if it is a complete
graph with loops at each node, and every edgeweight is positive. An
{\it unweighted graph} is a weighted graph where all the node- and
edgeweights are 1. Note that $e_G(S,T)$ reduces to the number of
edges in $G$ with one endnode in $S$ and the other in $T$ if $G$ is
unweighted.

Let $G$ be a graph and $k\ge 1$. The {\it $k$-fold blow-up} of $G$ is
the graph $G[k]$ obtained from $G$ by replacing each node by $k$
independent nodes, and connecting two new nodes if and only if their
originals were connected.  If $G$ is weighted, we define $G[k]$ to be
the graph on $nk$ nodes labeled by pairs $iu$, $i\in V(G)$,
$u=1,\dots, k$, with edgeweights $\beta_{iu,jv}(G[k])=\beta_{ij}(G)$
and nodeweights $\alpha_{iu}(G[k])=\alpha_i(G)$.  A related notion is
the notion of {\it splitting nodes}.  Here a node $i$ with nodeweight
$\alpha_i$ is replaced by $k$ nodes $i_1,\dots, i_k$ with nodeweights
$\alpha_{i_1},\dots,\alpha_{i_k}$ adding up to $\alpha_i$, with new
edgeweights $\beta_{i_u,j_v}=\beta_{i,j}$.  Up to a global rescaling of all
nodeweights, blowing up a graph by a factor $k$ is thus the same as
splitting all its nodes evenly into $k$ nodes, so that the new
weights $\alpha_{i_t}$ are equal to the old weights $\alpha_i$
divided by $k$.

As usual, a function from the set of simple graphs into the reals is
called a {\it simple graph parameter} if it is invariant under
relabeling of the nodes. Finally, we write $G\cong G'$ if $G$ and
$G'$ are isomorphic, i.e., if $G'$ can be obtained from $G$ by a
relabeling of its nodes.

\subsection{Homomorphism Numbers and Left Convergence}

Let $F$ and $G$ be two simple graphs. We  define $\hom(F,G)$ as
the number of homomorphisms from $F$ to $G$, i.e., the number of
adjacency preserving maps $V(F)\to V(G)$, and the {\it homomorphism
density} of $F$ in $G$ as
\[
t(F,G)=\frac {1}{|V(G)|^{|V(F)|}}\hom(F,G).
\]
The homomorphism density $t(F,G)$ is thus the probability that a
random map from $V(F)$ to $V(G)$ is a homomorphism.

Alternatively, one might want to consider the probability
$t_\inj(F,G)$ that a random injective map from $V(F)$ to $V(G)$ is
adjacency preserving, or the probability $t_\ind(F,G)$ that such a
map leads to an induced subgraph. Since most maps into a large graph
$G$ are injective, there is not much of a difference between $t(F,G)$
and $t_\inj(F,G)$. As for $t_\inj(\cdot,G)$ and $t_\ind(\cdot,G)$,
they can be quite different even for large graphs $G$, but by
inclusion-exclusion, the information contained in the two is strictly
equivalent. We therefore incur no loss of generality if we restrict
ourselves to the densities $t(\cdot,G)$.

We extend the  notion of homomorphism numbers to weighted graphs $G$
by setting
\begin{equation}
\label{hom-def} \hom(F,G)=\sum_{\phi:V(F)\to V(G)} \prod_{i\in V(F)}
\alpha_{\phi(i)}(G) \prod_{ij\in E(F)}\beta_{\phi(i),\phi(j)}(G)
\end{equation}
where the sum runs over all maps from $V(F)$ to $V(G)$,
and define
\begin{equation}
\label{t-def} t(F,G)=\frac{\hom(F,G)}{\alpha_G^k},
\end{equation}
where $k$ is the number of nodes in $F$.

It seems natural to think of two graphs $G$ and $G'$ as similar if
they have similar homomorphism densities.  This leads to the
following definition.

\begin{definition}
Let $(G_n)$ be a sequence of weighted graphs with uniformly bounded
edgeweights.  We say that $(G_n)$ is {\it convergent from the left},
or simply {\it convergent}, if $t(F,G_n)$ converges for any simple
graph $F$.
\end{definition}

In \cite{BCLSSV}, the definition of convergence was restricted to
sequences of graphs $(G_n)$ with $|V(G_n)|\to\infty$. As long as we
deal with simple graphs, this is reasonable, since there are only a
finite number of graphs with bounded size. But in this paper we also
want to cover sequences of weighted graphs on, say, the same set of
nodes, but with the node- or edgeweights converging; so we don't
assume that the number of nodes in convergent graph sequences tends
to infinity.

A simple example of a convergent graph sequence is a sequence of
random graphs $(G_{n,p})$, for which $t(F,G_{n,p})$ is convergent
with probability one, with $t(F,G_{n,p})\to p^{|E(F)|}$ as
$n\to\infty$. Other examples are quasirandom graph sequences $(G_n)$
for which $t(F,G_{n,p})\to p^{|E(F)|}$ by definition, and the
sequence of half-graphs $(H_{n,n})$, with $H_{n,n}$ defined as the
bipartite graph on $[2n]$ with an edge between $i$ and $j$ if $j\geq
n+i$ (see Examples \ref{RAND}, \ref{Q-RAND} and \ref{HALF} in
Section~\ref{sec:results-conv}).

\subsection{Cut-Distance}
\label{sec:results-Metric-Conv}

We define a notion of distance between two graphs, which will play a
central role throughout this paper. Among other
equivalences, convergence from
the left will be equivalent to convergence in this metric.

To illuminate the rather technical definition, we first define a
distance of two graphs $G$ and $G'$ in a special case, and then
extend it in two further steps. In all these definitions, it is not
hard to verify that the triangle inequality is satisfied.

\smallskip

(1) {\it $G$ and $G'$ are labeled graphs with the same set of unweighted nodes
$V$}. Several notions of distance appear in the literature, but for
our purpose the most useful is the cut or rectangle distance
introduced by Frieze and Kannan \cite{FK}:
\begin{equation}
\label{d-cut-unw-def} d_\square(G,G')=\max_{S,T\subset V} \frac
1{|V|^2}\Bigl|e_G(S,T)-e_{G'}(S,T)\Bigr|.
\end{equation}
The cut distance between two labeled graphs thus measures how
different two graphs are when considering the size of various cuts.
This definition can easily be generalized to weighted graphs $G$ and
$G'$ on the same set $V$, {\it and with the same nodeweights}
$\alpha_i=\alpha_i(G)=\alpha_i(G')$:
\begin{equation}\label{d-cut-w-def}
d_\square(G,G')=\max_{S,T\subset V}\frac1{\alpha_G^2}
\Bigl|e_G(S,T)-e_{G'}(S,T)\Bigr|.
\end{equation}

As a motivation of this notion, consider two independent random
graphs on $n$ nodes with edge density $1/2$. If we measure their
distance, say, by the number of edges we need to change to get one
from the other (edit distance), then their distance is very large
(with large probability). But the theory of random graphs
teaches us that these two graphs are virtually indistinguishable,
which is reflected by the fact that their $d_\square$ distance is
only $O(1/n)$ with large probability.

There are many other ways of defining or approximating the cut
distance; see Section \ref{MORE-NORMS}.

\smallskip

(2) {\it $G$ and $G'$ are unlabeled graphs with the same number of
unweighted nodes.} The cut-metric $d_\square(G,G')$ is not invariant
under relabeling of the nodes of $G$ and $G'$. For graphs without
nodeweights, this can easily be cured by defining a distance
$\iv\delta_\square(G,G')$ as the minimum over all ``overlays'' of
$G$ and $G'$, i.e.,
\begin{equation}
\label{delta-hat}
\iv\delta_\square(G,G')=\min_{\widetilde G \cong G}d_\square(\widetilde
G,G').
\end{equation}

\smallskip

(3) {\it $G$ and $G'$ are unlabeled graphs with different number of
nodes or with weighted nodes.} The distance notion \eqref{delta-hat}
does not extend in a natural way to graphs with nodeweights, since
it would not make much sense to overlay nodes with different
nodeweights, and even less
sense for graphs with different number of nodes.

To motivate the definition that follows, consider two graphs $G$ and
$G'$, where $G$ has three nodes with nodeweights equal to $1/3$, and
$G'$ has two nodes and nodeweights $1/3$ and $2/3$. Here a natural
procedure would be the following: first split the second node into
two nodes of weight $1/3$ and then calculated the optimal overlay of
the resulting two graphs on three nodes.

This idea naturally leads to the following notion of ``fractional
overlays'' of two weighted graphs $G$ and $G'$ on $n$ and $n'$ nodes,
respectively. Let us first assume that both $G$ and $G'$ have total
nodeweight $1$. Viewing $\alpha(G)$ and $\alpha(G')$ as probability
distributions, we then define a fractional overlay to be a coupling
between these two distributions. More explicitly, a {\it fractional
overlay} of $G$ and $G'$ is defined to be a nonnegative $n\times n'$
matrix $X$ such that
\[
\sum_{u=1}^{n'} X_{iu}=\alpha_i(G) \quad\text{and}\quad \sum_{i=1}^n
X_{iu}=\alpha_u(G').
\]
We denote the set of all fractional overlays by
$\XX(G,G')$.

Let $X\in\XX(G,G')$. Thinking of $X_{iu}$ as the portion of node $i$
that is mapped onto node $u$, we introduce the following ``
overlaid
graphs'' $G[X]$ and $G'[X\T]$  on $[n]\times[ n']$: in both $G[X]$
and $G'[X\T]$, the weight of a node $(i,u)\in [n]\times[ n']$ is
$X_{iu}$; in $G[X]$, the weight of an edge $((i,u),(j,v))$ is
$\beta_{ij}$, and in  $G'[X\T]$, the weight of an edge
$((i,u),(j,v))$ is $\beta'_{uv}$. Since $G[X]$ and $G'[X\T]$ have the
same nodeset, the distance $d_\square(G[X], G'[X\T])$ is now well
defined. Taking the minimum over all fractional overlays, this gives:

\begin{definition}
For two weighted graphs $G,G'$ with total nodeweight
$\alpha_G=\alpha_{G'}=1$, we set
\begin{equation}
\label{delta-altdef} \delta_\square(G,G')=\min_{X\in \XX(G,G')}
d_\square(G[X], G'[X\T]).
\end{equation}
If the total nodeweight of $G$ or $G'$ is different from $1$, we
define the distance between $G$ and $G'$ by the above formulas,
applied to the graphs $\tilde G$ and $\tilde G'$ obtained from $G$
and $G'$ by dividing all nodeweights by $\alpha_G$ and $\alpha_{G'}$,
respectively.
\end{definition}

Fractional overlays can be understood as integer overlays
of suitably blown up versions of $G$ and $G'$, at least if the
entries of $X$ are rational (otherwise, one has to take a limit of
blowups).  This observation shows that for two graphs $G$ and $G'$
with nodeweights one, we have
\begin{equation}
\label{delta=lim-hat-delta}
\delta_\square(G,G')
=\lim_{\bsa k,k'\to\infty\\ k/k'\,=\,n'/n\esa}
\iv\delta_\square(G[k],G'[k']).
\end{equation}

Note that $\delta_\square(G,G')$ can be $0$ for nonisomorphic graphs $G$
and $G'$; for example,
\begin{equation}\label{BLOW-DIST}
\delta_\square(G,G[k])=0
\end{equation}
for all $k\ge 1$. So $\delta_G$ is only a pre-metric; but we'll
call it, informally, a metric.

Of course, the definition of $\delta_\square(G,G')$ also applies if
$G$ and $G'$ have the same number of nodes, and it may give a value
different from $\iv\delta_\square(G,G')$. The following theorem
relates these two values.

\begin{theorem}\label{TWIN-SQUARE}
Let  $G_1$ and $G_2$ be two weighted graphs with edgeweights in
$[-1,1]$ and with the same number of unweighted nodes. Then
\[
\delta_\square(G_1,G_2)\le\iv\delta_\square(G_1,G_2)\le
{32}\delta_\square(G_1,G_2)^{{1/67}}.
\]
\end{theorem}

The first inequality is trivial, but the proof of the second is
quite involved, and will be given in Section \ref{sec:compdist}.

\subsection{Szemer\'edi Partitions of Graphs}
\label{sec:Graph-Szemeredi}

The Regularity Lemma of Szemer\'edi is a fundamental tool in graph
theory, which has a natural formulation in our framework,
as a result about approximating large graphs by small graphs.
Here we mostly use the so-called weak version due to Frieze and
Kannan \cite{FK}.

We need some notation: for a weighted graph $G$ and a partition
$\PP=\{V_1,\dots,V_k\}$ of $V(G)$, we define two weighted graphs
$G/\PP$ and $G_\PP$ as follows: Let $\alpha_{V_i}=\sum_{x\in
V_i}\alpha_x(G)$. The {\it quotient graph} $G/\PP$ is a weighted
graph on $[k]$, with nodeweights
$\alpha_i(G/\PP)=\alpha_{V_i}/\alpha_G$ and edgeweights
$\beta_{ij}(G/\PP) =\frac{e_G(V_i,V_j)}{\alpha_{V_i}\alpha_{V_j}}$,
while $G_\PP$ is a weighted graph on $V(G)$, with nodeweights
$\alpha_x(G_\PP)=\alpha_x(G)$ and edgeweights $\beta_{xy}(G_\PP)=
\beta_{ij}(G/\PP)$ for $x\in V_i$ and $y\in V_j$. These two graphs
have different number of nodes, but they are similar in the sense
that $\delta_\square(G/\PP,G_\PP)=0$.

In our language, the Weak Regularity Lemma of Frieze and Kannan
states that given a weighted graph $G$, one can find a partition
$\PP$ such that the graph $G_\PP$ is near to the original graph $G$
in the distance $d_\square$. We call the partition $\PP$ {\it weakly
$\eps$-regular} if $d_\square(G,G_\PP)\le\eps$. But for many
purposes, all that is needed is the quotient graph $G/\PP$. Since
$\delta_\square(G,G_\PP)\leq d_\square(G,G_\PP)$ and
$\delta_\square(G,G_\PP)=\delta_\square(G,G/\PP)$, the Weak
Regularity Lemma also guarantees a good approximation of the original
graph by a small weighted graph, the graph $H=G/\PP$. We summarize
these facts in the following lemma, which is essentially a
reformulation of the Weak Regularity Lemma of \cite{FK} in the
language developed in this paper.

\begin{lemma}[Weak Regularity Lemma \cite{FK}]\label{lem:szem-weak}
For every $\eps>0$, every weighted graph $G$ has a partition $\PP$ into at
most $4^{1/\eps^2}$ classes such that
\begin{equation}\label{szem-weak-d}
d_\square(G,G_\PP)\leq \eps\|G\|_2,
\end{equation}
so, in particular,
\begin{equation}
\label{szem-weak-delta} \delta_\square(G,G/\PP)\leq\eps\|G\|_2.
\end{equation}
\end{lemma}

In Lemma \ref{lem:szem-weak} we approximate the graph $G$ by a small
weighted graph $H=G/\PP$. If $G$ is simple, or more generally, has
edgeweights in $[0,1]$, it is possible to strengthen this by
requiring simple graphs $H$. Indeed, starting from a standard
strengthening of the Weak Regularity Lemma
(Corollary~\ref{cor:unif-Sem} (ii) below) to obtain a weighted graph
$H$ with nodeweights one, and then applying a simple randomization
procedure (Lemma \ref{lem:GH-CLOSE} below) to the edges of $H$ to
convert this graph into a simple graph, one gets the following lemma,
see appendix for details.

\begin{lemma}\label{lem:short-Szemeredi-simple}
Let $\eps>0$, let $q\geq 2^{20/\eps^2}$, and let $G$ be a weighted
graph with edgeweights in $[0,1]$. Then there exists a simple graph
$H$ on $q$ nodes such that $\delta_\square(G,H)\leq\eps$.
\end{lemma}

\subsection{Main Results}

\subsubsection{Left Convergence Versus Convergence in Metric}

Here we state one of the main results of this paper, namely, that
convergence from the left is equivalent to convergence in the metric
$\delta_\square$.

\begin{theorem}
\label{thm:left-equiv-metric} Let $(G_n)$ be a sequence of weighted
graphs with uniformly bounded edgeweights. Then $(G_n)$ is left
convergent if and only if it is a Cauchy sequence in the metric
$\delta_\square$.
\end{theorem}

In fact, we have the following quantitative version. To simplify our
notation, we only give this quantitative version for graphs with
edgeweights in $[-1,1]$; the general case follows by simply scaling
all edgeweights appropriately.

\begin{theorem}\label{F-CONV-CHAR}
Let $G_1,G_2$ be weighted graphs
with edgeweights in $[-1,1]$.

\smallskip

{\rm (a)}~Let $F$ be a simple graph, then
\[
|t(F,G_1)-t(F,G_2)|\le 4|E(F)|\delta_\square(G_1,G_2).
\]

\smallskip

{\rm (b)}~Let $k\geq 1$, and assume that
$
|t(F,G_1)-t(F,G_2)|\leq 3^{-{k^2}}
$
for every simple graph $F$ on $k$ nodes. Then
\[
\delta_\square(G_1,G_2) \le \frac{22}{\sqrt{\log_2 k}}.
\]
\end{theorem}

The first part of this theorem is closely related to the ``Counting
Lemma'' in the theory of Szemer\'edi partitions.
Theorems~\ref{thm:left-equiv-metric} and \ref{F-CONV-CHAR} will
follow from the analogous facts for graphons,
Theorems~\ref{LEFTCONV=DELTACONV} and \ref{W-CONV-CHAR}, see
Section~\ref{sec:Left-Metric-W}.

\subsubsection{Szemer\'edi Partitions for Graph Sequences}

Convergent graph sequences can also be characterized by the fact
that (for any
fixed error) they have Szemer\'edi partitions
which become more and more similar.

\begin{theorem}
\label{thm:Szem-Uniform} Let $(G_n)$ be a sequence of weighted
graphs with nodeweights $1$, edgeweights in $[-1,1]$, and
$|V(G_n)|\to\infty$ as $n\to\infty$.  Then $(G_n)$ is
left-convergent if and only if for every $\eps>0$ we can find an
integer $q\leq 2^{10/\eps^2}$, and a sequence of partitions,
$(\PP_n)$ such that the following two conditions hold.

(i) If $|V(G_n)|\geq q$, then $\PP_n$ is a weakly $\eps$-regular
partition of $G_n$ into $q$ classes.

(ii) As $n\to\infty$, the quotient graphs $G_n/\PP_n$ converge to a
weighted graph $H_\eps$ on $q$ nodes.
\end{theorem}

Note that the graphs in (ii) have the same node set $[q]$, so their
convergence to $H_\eps$ means simply that corresponding nodeweights and
edgeweights converge.

Let $G_n$ be a convergent sequence of weighted graphs obeying the
assumptions of this theorem. For $n$ sufficiently large, the quotient
graphs $G_n/\PP_n$ are then near to both the original graph $G_n$ and
the graph $H_\eps$, implying that $\delta_\square(G_n,H_\eps)\leq
2\eps$ whenever $n$ is large enough. Since $G_n$ is convergent, this
implies by Theorem \ref{thm:left-equiv-metric} that the graphs
$H_\eps$ form a convergent sequence as $\eps\to0$.

The theorem can be generalized in several directions.  First, using
the results of Section~\ref{sec:Norm-Conv}, we can relax the
condition on the nodeweights of $G_n$ to require only that $G_n$ has
no dominant nodeweights in the sense that the maximum nodeweight of
$G_n$ divided by the total nodeweight $\alpha_{G_n}$ goes to zero.
Second, we can strengthen the statement, to obtain a sequence of
partitions which satisfy the stronger regularity conditions of the
original Szemer\'edi Regularity Lemma \cite{Szem}.
We leave the details to the interested
reader, who will easily see how to modify the proof in
Section~\ref{sec:ConvSzemPart} to obtain these generalizations.

\subsubsection{Sampling}

Our above versions of Szemer\'edi's lemma (Lemmas~\ref{lem:szem-weak}
and \ref{lem:short-Szemeredi-simple}) state that any graph $G$ can be
well approximated by a small graph $H$ in the $\delta_\square$
distance.
While the proofs are
constructive, it will be very useful to know that such a small graph
can obtained by straightforward sampling. For simplicity, we state
the results for graphs with edgeweights in $[{-1},1]$. For a
graph $G$ and positive integer $n$, let $\RandG(n,G)$ denote the
(random) induced subgraph $G[S]$, where $S$ is chosen uniformly from
all subsets of $V(G)$ of size $n$.

\begin{theorem}\label{thm:sample}
Let $G$ be a weighted graph with nodeweights $1$ and edgeweights in
$[-1,1]$. Let $k\leq |V(G)|$. Then
\begin{equation}
\label{sample-bd} \delta_\square(G,\RandG(k,G)) \leq
\frac{10}{\sqrt{\log_2 k}}.
\end{equation}
with probability at least $1-e^{-k^2/(2\log_2 k)}$.
\end{theorem}

In order to prove this theorem, we will need a theorem which allows
us to compare samples from two weighted graphs on the same set of
nodes. This extends a result by Alon, Fernandez de la Vega, Kannan
and Karpinski \cite{AFKK,AFKK-Stoc}; in particular,
our result
concerns two graphs at arbitrary distance $d_\square(G_1,G_2)$, and
also gives an improvement in the error bound.

\begin{theorem}\label{thm:dist-test}
Let $G_1$ and $G_2$ be weighted graphs on a common vertex set $V$,
with nodeweights one and edgeweights in $[-1,1]$. Let $k\leq |V|$.
If $S$ is chosen uniformly from all subsets of $V$ of size $k$, then
\begin{equation}
\label{ProbG[X]-G} \Bigl|
d_\square(G_1[S],G_2[S])-d_\square(G_1,G_2) \Bigr|\leq
\frac{20}{k^{1/4}}
\end{equation}
with probability at least $1-2e^{-\sqrt k/8}$.
\end{theorem}

\subsubsection{Testing}

The above theorem allows us to prove several results for testing
graph parameters and graph properties in a straightforward way.

In this paper, we only consider parameter testing. We may want to
determine some parameter of $G$.  For example,
what is the edge density? Or how
large is the density of the maximum cut? Of course, we'll not be
able to determine the exact value of the parameter; the best we can
hope for is that if we take a sufficiently large sample, we can find
the approximate value of the parameter with large probability.

\begin{definition}\label{def:TESTABLE}
A graph parameter $f$ is {\it testable} if for every $\eps>0$ there
is a positive integer $k$ such that if $G$ is a graph with at least
$k$ nodes, then from the random subgraph $\RandG(k,G)$ we can compute
an estimate $\widetilde f$ of $f$ such that
\[
\Pr(|f(G)-\widetilde f|>\eps)\le \eps.
\]
It is an easy observation that we can always use
$\widetilde{f}=f(\RandG(k,G))$.
\end{definition}

Testability is related to our framework through the following
observation:

\begin{prop}\label{thm:TESTABLE-PAR}
{\rm (a)} A simple graph parameter $f$ is testable if and only if
$f(G_n)$ converges for every convergent graph sequence $(G_n)$ with
$|V(G_n)|\to\infty$.

\smallskip

{\rm (b)} A sequence $(G_n)$ of simple graphs with
$|V(G_n)|\to\infty$ is convergent if and only if $f(G_n)$ converges
for every testable simple graph parameter $f$.
\end{prop}

Using the notions and results concerning graph distance and
convergence above, we can give several characterizations of testable
parameters, see Section~\ref{sec:testing}.

Property testing, mentioned in the introduction, is related
to parameter testing in many ways. For example, Fischer and Newman
\cite{FN} proved that the edit distance (see Section
\ref{sec:RandWeightG}) from the set of graphs exhibiting a testable
property $\PP$ is a testable parameter. See also \cite{AFNS} for a
characterization of testable graph properties.

\section{Graphons} \label{sec:graphon}

In \cite{LSz1}, Lov\'asz and Szegedy introduced graphons as limits of
left-convergent graph sequences. Here we will first study the space
of graphons in its own right, defining in particular a generalization
of the distance $\delta_\square$ to graphons, and state the analogue
of Theorem~\ref{F-CONV-CHAR} for graphons. The discussion of graphons
as limit objects of left-convergent graph sequences will be postponed
to the last subsection of this section.

\subsection{Homomorphism Densities}
\label{sec:results-conv}

Let $\WW$ denote the space of all bounded measurable functions
$W:~[0,1]^2\to\R$ that are symmetric, i.e., $W(x,y)=W(y,x)$ for all
$x,y\in[0,1]$. Let {$\WWz$} be the set of functions $W\in\WW$ with
$0\le W(x,y)\le 1$. More generally, for a bounded interval
$I\subset\R$, let $\WW_I$ be the set of all functions $W\in\WW$ with
$W(x,y)\in I$. Given a function $W\in\WW$, we can think of the
interval $[0,1]$ as the set of nodes, and of the value $W(x,y)$ as
the weight of the edge $xy$.  We call the functions in $\WW$ {\it
graphons}.

We call a partition $\PP$ of $[0,1]$ {\it measurable} if all the
partition classes are (Lebesgue) measurable. The partition is an {\it
equipartition} if all of its classes have the same Lebesgue measure.
A {\it step function} is
a function $W\in\WW$ for which there is a
partition of $[0,1]$ into a finite number of measurable sets
$V_1,\dots,V_k$ so that $W$ is constant on every product set
$V_i\times V_j$. We call the sets $V_i$ the {\it steps} of the step
function. Often, but not always, we consider step functions whose
steps are intervals; we call these {\it interval step functions}. If
all steps of a step function have the same measure $1/k$, we say that
it has {\it equal steps}.

Every graphon $W$ defines a simple graph parameter as follows
\cite{LSz1}: If $F$ is a simple graph with $V(F)=\{1,\dots,k\}$,
then let
\begin{equation}
\label{t-lim-def} t(F,W)=\int_{[0,1]^k} \prod_{ij\in E(F)}
W(x_i,x_j)\,dx.
\end{equation}

Every weighted graph $G$ with nodes labeled $1,\dots,n$ defines an
interval step function $W_G$ such that $t(F,G)= t(F,W_G)$: We
scale the nodeweights of $G$ to sum to $1$. Let
$I_1=[0,\alpha_1(G)]$, $I_2=(\alpha_1(G),\alpha_1(G)+\alpha_2(G)]$,
$\dots$, and $I_n=[\alpha_1(G)+\dots+\alpha_{n-1}(G),1]$.  We then
set
\[
W_G(x,y)=\beta_{v(x)v(y)}(G)
\]
where $v(x)=i$ whenever $x\in I_i$. (Informally, we consider the
adjacency matrix of $G$, and replace each entry $(i,j)$ by a square
of size $\alpha_i\times\alpha_j$ with the constant function
$\beta_{ij}$ on this square.) If $G$ is unweighted, then the
corresponding interval step function is a $0$-$1$-function with
equal steps.

With this definition, we clearly have $W_G\in\WW$, and
\begin{equation}
\label{tG=tWG}
t(F,G)= t(F,W_G)
\end{equation}
for every finite graph $F$. The definition~\eqref{t-lim-def}
therefore gives a natural generalization of the homomorphism
densities defined in \eqref{hom-def} and \eqref{t-def}.

Using this notation, we can state the main result of
\cite{LSz1}:

\begin{theorem}\label{LSZTHM}
For every left-convergent sequence $(G_n)$ of simple graphs
there is a graphon $W$ with values in $[0,1]$ such that
\[
t(F,G_n)\to t(F,W)
\]
for every simple graph $F$. Moreover, for every graphon
$W$ with values in $[0,1]$ there is a left-convergent sequence of
graphs satisfying this relation.
\end{theorem}

\subsection{The Cut Norm for Graphons}
\label{sec:func-norm}

The distance $d_\square$ of two graphs introduced in
Section~\ref{sec:results-Metric-Conv} was extended to graphons by
Frieze and Kannan \cite{FK}. It will be given in terms of a norm on
the space $\WW$, the {\it rectangle} or {\it cut norm}
\begin{equation}\label{CUTNORM-W}
\begin{aligned}
\|W\|_\square &=\sup_{S,T\subseteq [0,1]} \biggl|\int_{S\times T}
W(x,y)\,dxdy\biggr| =\sup_{f,g:~[0,1]\to [0,1]}\biggl|\int
W(x,y)f(x)g(y)\,dx\,dy\biggr|
\end{aligned}
\end{equation}
where the suprema go over all pairs of measurable subsets and
functions, respectively. The cut norm is closely related to
$L_\infty\to L_1$ norm of $W$, considered as an operator on
$L^2([0,1])$:
\begin{align}\label{IONORM-W}
\|W\|_\IO &=\sup_{f,g:~[0,1]\to [-1,1]} \int_{[0,1]^2}
W(x,y)f(x)g(y)\,dx\,dy.
\end{align}
Indeed, the two norms are equivalent:
\begin{equation}\label{CUT-IO-W}
\frac{1}{4}\|W\|_\IO \le\|W\|_\square \le \|W\|_\IO.
\end{equation}
See section \ref{MORE-NORMS} for more on connections between the
cut norm and other norms.

It is  not hard to see that for any two weighted
graphs $G$ and $G'$
on the same set of nodes and with the same nodeweights,
\begin{equation}
\label{dG=dW_G}
d_\square(G,G')=\|W_G-W_{G'}\|_\square,
\end{equation}
where $W_G$ denotes the step function introduced in
Section~\ref{sec:results-conv}. The cut-norm therefore extends the
distance $d_\square$ from weighted graphs to graphons.

We will also need the usual norms of $W$ as a function from
$[0,1]^2\to\R$; we denote the corresponding $L_1$ and $L_2$ norms
(with respect to the Lebesgue measure) by $\|W\|_1$ and $\|W\|_2$.
The norm $\|.\|_2$ defines a Hilbert space, with inner product
\[
\langle U,W\rangle = \int_{[0,1]^2} U(x,y)W(y,x)\,dx\,dy.
\]

\subsection{Approximation by Step Functions}

We need to extend two averaging operations from graphs to graphons.
For $W\in\WW$ and every partition $\PP=(V_1,\dots,V_q)$ of $[0,1]$
into measurable sets, we define a weighted graph on $q$ nodes,
denoted by $W/\PP$, and called the quotient of $W$ and $\PP$,
by
setting $\alpha_i(W/\PP)=\lambda(V_i)$ (where $\lambda$ denotes the
Lebesgue measure) and
\begin{equation}\label{U-FACT}
\beta_{ij}(W/\PP)=\frac{1}{\lambda(V_i)\lambda(V_j)} \int_{V_i\times V_j}
W(x,y)\,dx\,dy.
\end{equation}
In addition to the quotient  $W/\PP$, we also consider the graphon
$W_\PP$ defined by
\begin{equation}\label{U-AVER}
W_\PP(x,y)= \beta_{ij}(W/\PP) \quad\text{whenever $x\in V_i$ and
$y\in V_j$.}
\end{equation}

It is not hard to check that the averaging operation $W\mapsto W_\PP$
is contractive with respect to the  norms introduced above:
\begin{equation}
\label{CONTRACT} \|W_\PP\|_\square\le \|W\|_\square, \quad
\|W_\PP\|_1\le \|W\|_1 \quad\text{and}\quad \|W_\PP\|_2\le\|W\|_2.
\end{equation}

The graphon $W_\PP$ is an approximation of $W$ by a step function with
steps $\PP$.  Indeed, it is the best such approximation, at least in
the $L_2$-norm:
\begin{equation}
\label{UP-Opt} \|W-W_\PP\|_2=\min_{U_\PP}\|W-U_\PP\|_2
\end{equation}
where the minimum runs over all step functions with steps $\PP$ (this
bound can easily be verified by varying the height of the steps in
$U_\PP$). While it is not true that $W_\PP$ is the best approximation
of $W$ with step in $\PP$ in the cut-norm, it is not off by more than
a factor of two, as observed in \cite{AFKK}:
\begin{equation}
\label{WP-2Opt} \|W-W_\PP\|_\square\leq 2
\min_{U_\PP}\|W-U_\PP\|_\square
\end{equation}
where the minimum runs over all step functions with steps $\PP$.
Indeed, combining the triangle inequality with the second bound in
\eqref{CONTRACT} and the fact that $(U_\PP)_\PP=U_\PP$, we conclude
that $\|W-W_\PP\|_\square\leq \|W-U_\PP\|_\square +
\|U_\PP-W_\PP\|_\square \leq \|W-U_\PP\|_\square +
\|U_\PP-W\|_\square$, as required.

The definition of $W_\PP$ raises the question on how well $W_\PP$
approximates $W$. One answer to this question is provided by
following lemma, which shows that $W$ can be approximated arbitrarily
well (pointwise almost everywhere) by interval step functions
with equal steps. (The lemma is an immediate consequence of the
almost everywhere differentiability of the integral function, see
e.g. Theorem 7.10 of \cite{Rud}.)

\begin{lemma}
\label{lem:as-approx} For a positive integer $n$, let $\PP_n$ be the
partition of $[0,1]$ into consecutive intervals of length $1/n$. For
any $W\in\WW$, we have $W_{\PP_n}\to W$ almost everywhere.
\end{lemma}

While the previous lemma gives the strong notion of almost everywhere
convergence, it does not give any bounds on the rate of convergence.
In this respect, the following lemma extending the weak Regularity
Lemma from graphs to graphons, is much better. In particular, it
gives a bound on the convergence rate which is independent of the
graphon $W$.

\begin{lemma}[\cite{FK}]\label{lem:W-szem}
For every graphon $W$ and every $\eps>0$, there exists a partition
$\PP$ of $[0,1]$ into measurable sets with at most
$4^{1/\eps^2}$
classes such that
\[
\|W-W_\PP\|_\square\le\eps\|W\|_2.
\]
\end{lemma}

With a slightly weaker bound for the number of classes, the lemma
follows from Theorem 12 of \cite{FK} and the bound \eqref{WP-2Opt},
or the results of \cite{LSz3}. As stated, it follows from
Lemma~\ref{lem:SZEM-ALG} in Section~\ref{sec:WRL-common}, which
generalizes both Lemma~\ref{lem:W-szem} and the analogous statement
for graphs, Lemma~\ref{lem:szem-weak}. The
Szemer\'edi Regularity Lemma \cite{Szem} also extends to
graphons in a straightforward way. See \cite{LSz3} for this and
further extensions.

At the cost of increasing the  bound on the number of classes,
Lemmas~\ref{lem:szem-weak} and \ref{lem:W-szem} can be strengthened
in several directions. In this paper, we need the following form,
which immediately follows from Lemmas~\ref{lem:szem-weak} and
\ref{lem:W-szem} by standard arguments,
see appendix for the details.

\begin{corollary}\label{cor:unif-Sem}
Let $\eps>0$, and $q\geq  2^{20/\eps^2}$. Then the following holds:

\smallskip

{\rm (i)} For all graphons $W$, there is an equipartition $\PP$ of
$[0,1]$ into $q$
measurable sets such that
\[
\|W-W_\PP\|_\square\leq \eps\|W\|_2.
\]
If we impose the additional constraint that $\PP$ refines a given
equipartition $\tilde\PP$ of $[0,1]$ into $k$ measurable sets,
it is possible to achieve this  bound provided $q$ is an integer
multiple of $k$ and $q/k\geq 2^{20/\eps^2}$.

\smallskip

{\rm (ii)} For all weighted graphs $G$ on at least $q$ nodes there
exists a partition $\PP=(V_1,\dots,V_q)$ of $V(G)$ such that
\[
d_\square(G,G_\PP)\leq \eps\|G\|_{2}
\]
and
\begin{equation}
\label{weight-cond} \Bigl|\sum_{u\in V_i}\alpha_u(G)-\frac
{\alpha_G}q\Bigr| <\alpha_{\max}(G) \qquad\text{for all
$i=1,\dots,q$}.
\end{equation}
\end{corollary}

\subsection{The Metric Space of Graphons}
\label{sec:FunctDist}

We now generalize the definition of the cut-distance
\eqref{delta-altdef}  from graphs to graphons.

Let $\MM$ denote the set of couplings of the uniform distribution on
$[0,1]$ with itself, i.e.,  the set of probability measures on
$[0,1]^2$ for which both marginals are the Lebesgue measure. (This
is the natural generalization of overlays from graphs to graphons).
We then define:
\[
\delta_\square(W,W')=\inf_{\mu\in\MM}\sup_{S,T\subseteq[0,1]^2}
\Bigl|\int\limits_{\ontop{(x,u)\in S}{(y,v)\in T}}
\bigl(W(x,y)-W'(u,v)\bigr)\,d\mu(x,u)\,d\mu(y,v)\Bigr|.
\]

For two step functions, finding the optimal ``overlay'' can be
described by specifying what fraction of each step of one function
goes onto each step of the other function, which amounts to
fractional overlay of the corresponding graphs. Hence the distances
of two unlabeled weighted graphs and the corresponding interval
step functions are the same:
\begin{equation}\label{STEPSAME}
\delta_\square(G,G') = \delta_\square(W_G,W_{G'}).
\end{equation}
It will be convenient to use the hybrid notation
$\delta_\square(U,G) = \delta_\square(U,W_G)$.

The next lemma gives an alternate representation of the distance
$\delta_\square(W,W')$, in which ``overlay'' is interpreted in terms
of measure-preserving maps rather than couplings. We need some
definitions. Recall that a map $\phi:~[0,1]\to[0,1]$ is {\it
measure-preserving}, if {the pre-image} $\phi^{-1}(X)$ is measurable
for every measurable set $X$, and
$\lambda(\phi^{-1}(X))=\lambda(X)$. A {\it measure-preserving
bijection} is a measure-preserving map whose inverse map exists and
is also measurable (and then also measure-preserving). Finally, we
consider certain very special measure-preserving maps defined as
follows: Let us consider the partition $\PP_n=(V_1,\dots,V_n)$ of
$[0,1]$ into consecutive intervals of length $1/n$, and let $\pi$ be
a permutation of $[n]$. Let us map each $V_i$ onto $V_{\pi(i)}$ by
translation, to obtain a piecewise linear measure-preserving map
$\tilde{\pi}:~[0,1]\to[0,1]$. We call $\tilde{\pi}$ an {\it $n$-step
interval permutation}.

For $W\in \WW$ and $\phi:~[0,1]\to[0,1]$, we define $W^\phi$ by
$W^\phi(x,y)=W(\phi(x),\phi(y))$.

\begin{lemma}
\label{lem:delta=delta} Let $U,W\in \WW$. Then
\begin{align}
\delta_\square(U,W)
&=\inf_{\phi,\psi}\|U^\phi-W^\psi\|_\square\label{eq:DD1}\\
&(\text{where the infimum is over all measure-preserving maps
$\phi,\psi:~[0,1]\to[0,1]$})\nonumber\\
&=\inf_{\psi}\|U-W^\psi\|_\square\label{eq:DD2}\\
&(\text{where the infimum is over all measure-preserving bijections
$\psi:~[0,1]\to[0,1]$})\nonumber\\
&=\lim_{n\to\infty}\min_{\pi}\|U-W^{\tilde{\pi}}\|_\square\label{eq:DD3}\\
&(\text{where the minimum is over all permutations $\pi$ of
$[n]$.})\nonumber
\end{align}
\end{lemma}

The proof of the lemma is somewhat tedious, but straightforward,
see appendix for details.

Note that, for $W\in\WW$, for an $n$-step interval permutation
$\tilde{\pi}$, and for the partition $\PP_n$ of $[0,1]$ into
consecutive intervals of lengths $1/n$, the graph
$W^{\tilde{\pi}}/\PP_n$ is obtained from $W/\PP_n$ by a permutation
of the nodes of $W/\PP_n$.  As a consequence, the identity
\eqref{eq:DD3} is equivalent to the following analogue of
\eqref{delta=lim-hat-delta} for graphons:
\begin{equation}
\delta_\square(U,W)=
\lim_{n\to\infty}\iv\delta_\square(U/\PP_n,W/\PP_n).
\end{equation}

Using Lemma~\ref{lem:delta=delta}, it is easy to verify that
$\delta_\square$ satisfies the triangle inequality. Strictly
speaking, the function $\delta_\square$ is just a pre-metric, and not
a metric: Formula \eqref{eq:DD2} implies that
$\delta_\square(W,W')=0$ whenever $W=W^\phi$ for some measure
preserving transformation $\phi:~[0,1]\to[0,1]$. Nevertheless, for
the sake of linguistic simplicity, we will often refer to it as a
distance or metric, taking the implicit identification of graphs or
graphons with distance zero for granted. Note that the fact that
there are graphons $W,W'\in\WW$ which are different but have
``distance'' zero is not just a peculiarity of the limit: for simple
graphs $\delta_\square(G,G')$ is zero if, e.g., $G'$ is a blow-up of
$G$ (cf. \eqref{BLOW-DIST}). We'll say more about graphons with
distance $0$ in the next section.

It was proved in \cite{LSz3} that the metric space
$(\WWz,\delta_\square)$ is compact. Since $\WW_I$ and $\WWz$
are linear images of each other, this immediately implies the
following proposition.

\begin{prop}\label{thm:compact}
Let $I$ be a finite interval, and let $\WW_I$ be the set of graphons
with values in $I$. After identifying graphons with $\delta_\square$
distance zero, the metric space $(\WW_I,\delta_\square)$ is compact.
\end{prop}

\subsection{Left Versus Metric Convergence}\label{sec:Left-Metric-W}

We are now ready to state the analogue of
Theorems~\ref{thm:left-equiv-metric} and \ref{F-CONV-CHAR}.
for graphons. We start with the analogue of
Theorem~\ref{F-CONV-CHAR}.

\begin{theorem}\label{W-CONV-CHAR}
Let $W,W'\in\WW$, let $C=\max\{1,\|W\|_\infty,\|W'\|_\infty\}$, and
let $k\geq 1$.

\smallskip

{\rm (a)}~If  $F$ is a simple graph with $m$ edges, then
\begin{equation}
\label{LEFTDIST-delta} |t(F,W)-t(F,W')|\le 4 mC^{m-1}{
\delta_\square(W,W')}.
\end{equation}

\smallskip

{\rm (b)}~If
$
|t(F,W)-t(F,W')|\leq 3^{-{k^2}}
$
for every simple graph $F$ on $k$ nodes, then
\[
\delta_\square(W, W') \le \frac{22C}{\sqrt{\log_2 k}}.
\]
\end{theorem}

The first statement of this theorem is closely related to the
``Counting Lemma'' in the theory of Szemer\'edi partitions, and
gives an extension of a similar result of \cite{LSz1} for functions
in $\WWz$ to general graphons. It shows that for any simple graph
$F$, the function $W\mapsto t(F,W)$ is Lipschitz-continuous in the
metric $\delta_\square$, and is reasonably easy to prove. By
contrast, the proof of the second one is more involved and relies on
our results on sampling.

In particular, we will need an analogue of Theorems~\ref{thm:sample}
and \ref{thm:dist-test} to sampling from graphons. These theorems
are stated and proved in Section~\ref{sec:dist-sample-W}
(Theorem~\ref{thm:sample2} and \ref{thm:NormSample}). Using
Theorem~\ref{thm:sample2}, we then prove Theorem~\ref{W-CONV-CHAR}
in  Section~\ref{sec:proof-main-W}.

Theorem~\ref{W-CONV-CHAR} immediately implies the analogue
of Theorem~\ref{thm:left-equiv-metric} for graphons:

\begin{theorem}\label{LEFTCONV=DELTACONV}
Let $I$ be a finite interval and let $(W_n)$ be a sequence of
graphons in $\WW_I$. Then the following are equivalent:

\smallskip

{\rm (a)} $t(F,W_n)$ converges for all finite simple graphs $F$;

\smallskip

{\rm (b)} $W_n$ is a Cauchy sequence in the $\delta_\square$ metric;

\smallskip

{\rm (c)} there exists a $W\in \WW_I$ such that $t(F,W_n)\to t(F,W)$
for all finite simple graphs $F$.

\smallskip

Furthermore,  $t(F,W_n)\to t(F,W)$ for all finite simple graphs $F$
for some $W\in \WW$ if and only if $\delta_\square(W_n,W)\to 0$.
\end{theorem}

Note that by equations \eqref{tG=tWG} and \eqref{STEPSAME},
Theorems~\ref{thm:left-equiv-metric} and \ref{F-CONV-CHAR}
immediately follow from Theorems~\ref{LEFTCONV=DELTACONV} and
\ref{W-CONV-CHAR}. Together with Theorem \ref{LSZTHM}, these
results imply that after identifying graphons with distance zero, the
set of graphons $\WWz$ is the completion of the metric space of
simple graphs. Proposition~\ref{thm:compact}, equation
\eqref{STEPSAME} and Lemma~\ref{lem:as-approx} easily imply that the
existence of the limit object (Theorem \ref{LSZTHM}) can be
extended to convergent sequences of weighted graphs:

\begin{corollary}\label{thm:LIMIT-EXIST}
For any convergent sequence $(G_n)$ of weighted graphs with
uniformly bounded edgeweights there exists a graphon $W$ such that
$\delta_\square(W_{G_n},W)\to 0$. Conversely, any graphon $W$ can be
obtained as the limit of a sequence of weighted graphs with
uniformly bounded edgeweights. The limit of a convergent graph
sequence is essentially unique: If $G_n\to W$, then also $G_n\to W'$
for precisely those graphons $W'$ for which
$\delta_\square(W,W')=0$.
\end{corollary}

As another consequence of Theorem \ref{LEFTCONV=DELTACONV}, we
get a characterization of graphons of distance $0$:

\begin{corollary}\label{cor:D=0}
For two graphons $W$ and $W'$ we have $\delta_\square(W,W')=0$ if
and only if $t(F,W)=t(F,W')$ for every simple graph $F$.
\end{corollary}

Another characterization of such pairs is given in \cite{BCL}:
$\delta_\square(W,W')=0$ if and only if there exists a third graphon
$U$ such that $W=U^\phi$ and $W'=U^\psi$ for two measure-preserving
functions $\phi,\psi:~[0,1]\to[0,1]$ (in other words, the infimum in
\eqref{eq:DD1} is a minimum if the distance is $0$).

\subsection{Examples}

\begin{example}[Random graphs]\label{RAND}
Let $\bG(n,p)$ be a random graph on $n$ nodes with edge density $0\le
p\le 1$; then it is not hard to prove (using high concentration
results) that the sequence $(\bG(n,p), n=1,2,\dots)$ is convergent
with probability 1. In fact, $t(F,\bG(n,p){)}$ converges to
$p^{|E(F)|}$ with probability 1, and so (with probability $1$)
$\bG(n,p)$ converges to the constant function $W=p$.
\end{example}

\begin{example}[Quasirandom graphs]\label{Q-RAND}
A graph sequence is quasirandom with density $p$ if and only if it
converges to the constant function $p$.  Quasirandom graph
sequences have many other interesting characterizations in terms of
edge densities of cuts, subgraphs, etc. \cite{CGW}. In the second
part of this paper we'll discuss how most of these characterizations
extend to convergent graphs sequences.
\end{example}

\begin{example}[Half-graphs]\label{HALF}
Let $H_{n,n}$ denote the bipartite graph on $2n$ nodes
$\{1,\dots,n,1',\dots,n'\}$, where $i$ is connected to $j'$ if and
only if $i\le j$. It is easy to see that this sequence is convergent,
and its limit is the function
\[
W(x,y)=
  \begin{cases}
    1, & \text{if $|x-y|\ge 1/2$}, \\
    0, & \text{otherwise}.
  \end{cases}
\]
\end{example}

\begin{example}[Uniform attachment]\label{U-GROW}
Various sequences of growing graphs, motivated by (but different
from) internet models, are also convergent. We define a
{(dense)}
{\it  uniform attachment graph sequence} as follows:
if we have a current graph $G_n$ with $n$
nodes, then we create a new isolated node, and then for every pair of
previously nonadjacent nodes, we connect them with probability $1/n$.

One can prove that with probability $1$, the sequence $(G_n)$ has a
limit, which is the function $W(x,y)=\min(x,y)$. From this, it is
easy to calculate that with probability $1$, the edge density of
$G_n$ tends to $\int W=1/3$. More generally, the density of copies of
any fixed graph $F$ in $G(n)$ tends (with probability $1$) to
$t(F,W)$, which can be evaluated by a simple integration.
\end{example}

\section{Sampling}
\label{sec:sampling}

\subsection{Injective and Induced Homomorphisms}

In order to discuss sampling, we will consider not only the number
of homomorphisms defined earlier, but also the number of injective
and induced homomorphisms between two simple graphs $F$ and $G$. We
use $\inj(F,G)$  to denote the number of {\it injective}
homomorphisms from $F$ to $G$, and  $\ind(F,G)$ to denote the number
of those injective homomorphisms that also preserve non-adjacency
(equivalently, the number of embeddings of $F$ into $G$ as an
induced subgraph).

{We will need to generalize these notions to the case where $G$
is a weighted graph with nodeweights one and edgeweights
$\beta_{ij}(G)\in\R$, where we define
\begin{equation}
\label{inj-def} \inj(F,G)=\sum_{\phi\in\Inj(F,G)} \prod_{uv\in E(F)}
\beta_{\phi(u),\phi(v)}(G)
\end{equation}
and
\begin{equation}
\label{ind-def} \ind(F,G)= \sum_{\phi\in\Inj(F,G)} \prod_{uv\in E(F)}
\beta_{\phi(u),\phi(v)}(G) \prod_{uv\in E(\overline{F})}
\Bigl(1-\beta_{\phi(u),\phi(v)}(G)\Bigr).
\end{equation}
Here $\Inj(F,G)$ denotes the set of injective maps from $V(F)$ to
$V(G)$, and $E(\overline{F})$ consists of all pairs $\{u,v\}$ of
distinct nodes such that $uv\notin E(F)$. We also introduce the
densities
\begin{equation}
t_\inj(F,G)=\frac{\inj(F,G)}{(|V(G)|)_{|V(F)|}}
\quad\text{and}\quad
t_\ind(F,G)=\frac{\ind(F,G)}{(|V(G)|)_{|V(F)|}}
\end{equation}
where $(n)_k=n(n-1)\cdots(n-k+1)$.
}

The quantities $t(F,G),t_\inj(F,G)$ and $t_\ind(F,G)$ are closely
related. It is easy to see  that
\begin{equation}\label{T0T1}
t_\inj(F,G)=\sum_{F'\supset F} t_{\ind}(F',G)
\qquad\text{and}\qquad
t_\ind(F,G)=\sum_{F'\supset F}(-1)^{|E(F')\setminus E(F)|}
t_\inj(F',G)
\end{equation}
whenever $F$ is simple and $G$ is a weighted graph with nodeweights
$\alpha_i(G)=1$. The quantity $t(F,G)$ is not expressible as a
function of the values $t_\inj(F,G)$ (or $t_\ind(F,G)$), but for
large graphs, they are essentially the same.  Indeed, bounding the
number of non-injective homomorphisms from $V(F)$ to $V(G)$ by
$\binom{|V(F)|}{2}|V(G)|^{|V(F)|-1}$, one easily proves that
\begin{equation}\label{TT0}
|t(F,G) -t_\inj(F,G)|
<\frac{2}{|V(G)|}\binom{|V(F)|}{2}\|G\|_\infty^{|E(F)|}.
\end{equation}
If all edgeweights of $G$ lie in the interval $[0,1]$, this bound can
be strengthened to
\begin{equation}\label{TT0-better}
|t(F,G) -t_\inj(F,G)| <\frac{1}{|V(G)|}\binom{|V(F)|}{2};
\end{equation}
see \cite{LSz1} for a proof.

\subsection{Sampling Concentration}

We will repeatedly use the following consequences of Azuma's
Inequality:

\begin{lemma}\label{lem:Azuma}
Let $(\Omega,\AA,\Prob)$ be a probability space, let
$k$ be a positive integer, and let $C>0$.

(i) Let $Z=(Z_1,\dots,Z_k)$, where $Z_1,\dots,Z_k$ are independent
random variables, and $Z_i$ takes values in some measure space
$(\Omega_i,\AA_i)$. Let $f:~\Omega_1\times\cdots\times\Omega_k\to \R$
be a measurable function. Suppose that $|f(x)-f(y)|\le C$ whenever
$x=(x_1,\dots,x_k)$ and $y=(y_1,\dots,y_k)$ differ only in one
coordinate.  Then
\begin{equation}
\Prob\Bigl(f(Z)>\E[f(Z)]+\lambda C\Bigr) <
e^{-\lambda^2/2k} \qquad\text{and}\qquad
\Prob\Bigl(\bigl|f(Z)-\E[f(Z)]\bigr|>\lambda C\Bigr) <
2e^{-\lambda^2/2k}.
\label{Azuma1}
\end{equation}

(ii) The bounds \eqref{Azuma1} also hold if $Z_1,\dots,Z_k$ are
chosen uniformly without replacement from some finite set $V$ and
$|f(x)-f(y)|\leq C$ for all $x$ and $y$ which either differ in at
most one component, or can be obtained from each other by permuting
two components.

\end{lemma}

\begin{proof}
Let $M_j=\E[\,f(Z)\mid Z_1,\dots,Z_j]$.  Then $M_0,\dots,M_k$ is a
martingale with bounded martingale differences (for the case (ii)
this requires a little calculation which we leave to the reader). The
statement now follows from Azuma's inequality for bounded
martingales.
\end{proof}

\subsection{Randomizing Weighted Graphs}
\label{sec:RandWeightG}

Given a weighted graph $H$ with nodeweights $1$ and edgeweights in
$[0,1]$, let $\bG(H)$ denote the random simple graph with $V(G)=V(H)$
obtained by joining nodes $i$ and $j$ with probability
$\beta_{ij}(H)$ (making an independent decision for every pair $ij$,
and ignoring the loops in $H$).

We need two simple properties of this well known construction. To
state the first, we define the edit distance $d_1$ of two weighted
graphs with the same node set $[n]$ and nodeweights $1$ as
\[
d_1(H_1,H_2)=\frac 1{n^2}\sum_{i,j=1}^n|\beta_{ij}(H_1)-\beta_{ij}(H_2)|.
\]

\begin{lemma}\label{lem:GH-COUPL}
Let $H_1$ and $H_2$ be two weighted graphs on the same set of nodes
with nodeweights $1$ and with edgeweights in $[0,1]$. Then
$\bG(H_1)$ and $\bG(H_2)$ can be coupled so that
\[
\E(d_1(\bG(H_1),\bG(H_2))) = d_1(H_1,H_2).
\]
\end{lemma}

\begin{proof}
For every edge $ij$, we couple the decisions about the edge so that
in both graphs this edge is inserted with probability
$\min(\beta_{ij}(H_1),\beta_{ij}(H_2))$ and missing with probability
$1-\max(\beta_{ij}(H_1),\beta_{ij}(H_2))$. So the probability that
the edge is present in exactly one of $\bG(H_1)$ and $\bG(H_2)$ is
$|\beta_{ij}(H_1)-\beta_{ij}(H_2)|$, which proves the lemma.
\end{proof}

\begin{lemma}\label{lem:GH-CLOSE}
Let $H$ be a weighted graph on $n$ nodes
with nodeweights $1$ and with edgeweights in $[0,1]$. Then
\[
\Pr\biggl(d_\square(H,\bG(H)) < \frac {4}{\sqrt n}
\biggr)
>1-2^{-n}.
\]
\end{lemma}

\begin{proof}
Let $V(H)=V(\bG(H))=V$, let $\mu=3/\sqrt n$, and let
$\widetilde H$ be the graph obtained from $H$ by deleting all
diagonal entries in $\beta(H)$. Fix two sets $S,T\subseteq V$. For
$i\neq j\in V$, let
\[
X_{ij}=
  \begin{cases}
    1 & \text{if $ij\in E(\bG(H))$}, \\
    0 & \text{otherwise}.
  \end{cases}
\]
Observe that the expectation of $e_{\bG(H)}(S,T)$ is $e_{\widetilde
H}(S,T)$. Since $e_{\bG(H)}(S,T)$ is a function of the
$\frac{n(n-1)}2$ independent random variables $(X_{ij})_{i< j}$ that
changes by at most $2$ if we change one of these variables, we may
apply Lemma~\ref{lem:Azuma} to conclude that
\[
\Pr\Bigl(\bigl|e_{\bG(H)}(S,T)-e_{\widetilde H}(S,T)\bigr| \geq\mu
n^2\Bigr)) \le 2\exp\Bigl(-\frac{{\mu}^2 n^4}{4n(n-1)}\Bigr) <
\exp\Bigl(-\frac{\mu^2 n^2}4\Bigr).
\]
(Here we used that $e^{-\mu^2 n/4}<1/2$ in the last step).
Taking into account that the number of pairs $(S,T)$ is $4^n$, we
concluded that the probability that $d_\square(\widetilde
H,\bG(H))<\mu=3/\sqrt n$
is larger than $1-4^ne^{-{\mu}^2 n^2/{4}}
\geq 1-2^{-n} $. Since
$d_\square(\widetilde H,H)\leq 1/n\leq 1/\sqrt n$, this
completes the proof.
\end{proof}

\subsection{$W$-random Graphs}

Given a graphon $W\in\WW$ and a subset $S\subseteq [0,1]$, we define
the weighted graph $W[S]$ on node set $S$, all nodes with weight $1$,
in which $\beta_{xy}(W[S])=W(x,y)$. If $W\in\WWz$, then we can
construct a random {\it simple} graph $\iv{W}[S]$ on $S$ by
connecting nodes $X_i$ and $X_j$ with probability $W(X_i,X_j)$
(making an independent decision for every pair).

This construction gives rise to two random graph models defined by
the graphon $W$. For every integer $n>0$, we generate a {\it
{\rm $W$-}random weighted graph} $\bH(n,W)$ on nodes $\{1,\dots,n\}$ as
follows: We generate $n$ independent samples $X_1,\dots,X_n$ from the
uniform distribution on $[0,1]$, and consider $W[\{X_1,\dots,X_n\}]$
(renaming $i$ the node $X_i$). If $W\in\WWz$, then we also define
the {\it {\rm $W$-}random (simple) graph}
$\bG(n,W)\cong\iv{W}[\{X_1,\dots,X_n\}]$.

{When proving concentration, it will often be useful to generate
$\bG(n,W)$ by first independently choosing $n$ random
variables $X_1,\dots,X_n$ and $n(n+1)/2$ random variables
$Y_{ij}$ ($i\leq j$) uniformly at random from $[0,1]$,
and then defining $\bG(n,W)$ to be the graph with an edge
between $i$ and $j$ whenever $Y_{ij}\leq W(X_i,X_j)$.
This allows us to express the adjacency matrix of
$\bG(n,W)$ as a function of the {\it independent}
random variables $Z_1=(X_1,Y_{11})$,
$Z_2=(X_2,Y_{12},Y_{22})$, \dots
$Z_n=(X_n,Y_{1n},Y_{2n},\dots,Y_{nn})$, as required for the application of
Lemma~\ref{lem:Azuma} (i).}

It is easy to see that for every simple graph $F$ with $k$ nodes
\begin{equation}
\label{W-random-Exp} \E\bigl(t_\inj(F,\bG(n,W))\bigr)=
\E\bigl(t_\inj(F,\bH(n,W))\bigr)= t(F,W),
\end{equation}
where the second equality holds for all $W$ while the first requires
$W\in\WWz$. From this we get that
\begin{equation}
\label{random-W-Exp-t}
\bigl|\E\bigl(t(F,\bH(n,W))\bigr)- t(F,W)\bigr|<\frac{
2}n\binom{k}{2}
\qquad\text{if}\qquad {\|W\|_\infty\leq 1,}
\end{equation}and
\begin{equation}
\label{W-random-Exp-t}
\bigl|\E\bigl(t(F,\bG(n,W))\bigr)-
t(F,W)\bigr|<\frac{1}n\binom{k}{2}
\qquad\text{if}\qquad{W\in\WWz.}
\end{equation}

Concentration for the $W$-random graph $\bG(n,W)$ was established in
Theorem 2.5 of \cite{LSz1}. To get concentration for $W$-weighted
random graphs, we use Lemma \ref{lem:Azuma}. This gives the following
lemma, which also slightly improves the bound of Theorem 2.5 of \cite{LSz1}
for $W$-random graphs.
See the appendix for details of the proof.

\begin{lemma}\label{lem:t-conc}
Let $F$ be a simple graph on $k$ nodes, let $0<\eps<1$
and let $W\in\WW$.  Then
\begin{equation}
\label{WH-Concentr} \Prob\Bigl(|t(F,\bH(n,W))-t(F,W)| > \eps
\Bigr) \le 2\exp\left(-\frac{\eps^2}{{
11}k^2}n\right)
\qquad\text{if}\qquad {\|W\|_\infty\leq 1,}
\end{equation}
and
\begin{equation}
\label{WG-Concentr} \Prob\Bigl(|t(F,\bG(n,W))-t(F,W)| > \eps \Bigr)
\le 2\exp\left(-\frac{\eps^2}{4k^2}n\right)
\qquad\text{if}\qquad {W\in\WWz.}
\end{equation}
\end{lemma}

From this lemma we immediately get the following:

\begin{theorem}\label{WRAND-CONV}
{\rm (a)} For any $W\in\WW$, the graph sequence $\bH(n,W)$ is
convergent with probability $1$, and its limit is the graphon $W$.

{\rm (b)} For any $W\in\WWz$, the graph sequence $\bG(n,W)$ is
convergent with probability $1$, and its limit is the graphon $W$.
\end{theorem}

\subsection{The Distance of Samples}
\label{sec:dist-sample-W}

The closeness of a sample to the original graph lies at the heart of
many results in this and the companion paper \cite{dense2}.  We will
prove these results starting with an extension of Theorem
\ref{thm:dist-test}, which is an improvement of a theorem of Alon,
Fernandez de la Vega, Kannan and Karpinski \cite{AFKK,AFKK-Stoc},
as discussed earlier.

\begin{theorem}\label{thm:NormSample}
Let $k$ be a positive integer.

\smallskip

{\rm (i)} If $U\in\WW$, then
\begin{equation}
\label{ProbH-UpBd}
\Bigl|\|\bH(k,U)\|_\square -\|U\|_\square\Bigr|
\leq \frac {10} {k^{1/4}}\|U\|_\infty
\end{equation}
with probability at least $1-2e^{-\sqrt k/8}$.

\smallskip

{\rm (ii)} If $U_1,U_2\in\WWz$, then $\bG(k,U_1)$ and $\bG(k,U_2)$
can be coupled in such a way that
\begin{equation}
\label{ProbG1G2-diff}
\Bigl|
d_\square(\bG(k,U_1),\bG(k,U_2))-\|U_1-U_2\|_\square
\Bigr|\leq \frac{10}{k^{1/4}},
\end{equation}
with probability at least $1-e^{-\sqrt k/8}$.
\end{theorem}

\begin{proof}
(i) Since the assertion is homogeneous in $U$, we may assume without
loss of generality that $\|U\|_\infty=1$. The proof proceeds in two
steps: first we prove that for $k\geq 10^4$, we have
\begin{equation}
\label{ExpH-UpBd} -\frac 2k \leq
\E\bigl(\|\bH(k,U)\|_\square\bigr)-\|U\|_\square < \frac 8{k^{1/4}},
\end{equation}
and then use Lemma~\ref{lem:Azuma} to prove concentration.

It turns out that the most difficult step is the proof of an upper
bound on the expectation of $\|\bH(k,U)\|_\square $.  To obtain this
bound, we use a refinement of the proof strategy of
\cite{AFKK}, using  Lemma 3 from that paper as our  starting
point. The main difference between their and our proofs is that we
first bound the expectation of $\|\bH(k,U)\|_\square $, and only  use
Lemma~\ref{lem:Azuma} in the very end. This allows us to simplify
their proof, giving at the same time a better dependence of our
errors on $k$.

Let $X_1,\dots,X_k$ be i.i.d.~random variables, distributed uniformly
in $[0,1]$, let $B=B({\bX})$ be the $k\times k$ matrix with entries
$B_{ij}=U(X_i,X_j)$, and for $S_1,S_2\subset [k]$, let $B(S_1,S_2)=
\sum_{i\in S_1,j\in S_2} B_{ij}$. Finally, given a set $S\subset
[k]$, let $P(S)$ be the set of nodes $i\in [k]$ such that
$B(\{i\},S)>0$, and $N(S)$ be the set of nodes for which
$B(\{i\},S)<0$.  We will prove upper and lower bounds on the
expectation of $||B\|_\square =\max_{S_1,S_2\subset
[k]}|B(S_1,S_2)|$.

The lower bound is a simple sampling argument: consider two
measurable subsets $T_1,T_2\subset [0,1]$, and let $S_1=\{i\in
[k]\colon X_i\in T_1\}$, and similarly for $S_2$.  Then
\[
\begin{aligned}
\E\|\bH(k,U)\|_\square &\geq\frac
1{k^{2}}\Bigl|\E\Bigl[B(S_1,S_2)\Bigr]\Bigr|
= \Bigl| \frac{k-1}k\int_{T_1\times T_2} U(x,y)\,dx\,dy + \frac
1k\int_{T_1\cap T_2} U(x,x)\,dx\Bigr|
\\
&\geq \Bigl| \int_{T_1\times T_2} U(x,y)\,dx\,dy \Bigr|-\frac 2k.
\end{aligned}
\]
Taking the supremum over all measurable sets
${S}_1,{S}_2\subset [0,1]$, this proves the lower bound in
\eqref{ExpH-UpBd}.

To prove an upper bound on the expectation of
$\|\bH(k,U)\|_\square$, we start from Lemma 3 of \cite{AFKK}.  In
our context, it states that for a random subset $Q\subset [k]$ of
size $p$, we have
\[
B(S_1,S_2)\leq \E_Q[B(P(Q\cap S_2),S_2)]
+\frac k{\sqrt p}\|B\|_F,
\]
where $\|B\|_F=\sqrt{\sum_{i,j}B_{ij}^2}\le k$. Inserting this
inequality into itself, we obtain that
\[
B(S_1,S_2) \leq \E_{Q_1,Q_2}\Bigl[ \max_{T_i\subseteq Q_i}
B(P(T_1),P(T_2))\Bigr] +\frac {2k^2}{\sqrt p}.
\]

In order to take the expectation over the random variables
$X_1,\dots,X_k$, it will be convenient to decompose $P(T_1)$ and
$P(T_2)$ into the parts which intersect $Q=Q_1\cup Q_2$ and the parts
which intersect  $Q^c=[k]\setminus (Q_1\cup Q_2)$. Let
$P_Q(T)=P(T)\cap Q$ and $P_{Q^c}(T)=P(T)\setminus Q$. Since
$|P(T)|\leq k$ and $|P(T)\setminus P_{Q^c}(T)|=|P_Q(T)|\leq 2p$ we
have that $B(P(T_1),P(T_2))\leq B(P_{Q^c}(T_1),P_{Q^c}(T_2))+4pk$,
implying that
\[
B(S_1,S_2) \leq \E_{Q_1,Q_2}\Bigl[ \max_{T_i\subseteq Q_i}
B(P_{Q^c}(T_1),P_{Q^c}(T_2))\Bigr] +\frac {2k^2}{\sqrt p} + 4pk.
\]
Applying this estimate to $-B$ and taking the maximum of the two
bounds, this gives
\begin{equation}
\label{B+bound} \|B\|_\square \leq \E_{Q_1,Q_2}\Bigl[
\max_{T_i\subseteq Q_i} \max\{B(P_{Q^c}(T_1),P_{Q^c}(T_2)),
-B(N_{Q^c}(T_1),N_{Q^c}(T_2))\} \Bigr] +\frac {2k^2}{\sqrt p} + 4pk,
\end{equation}
where $N_{Q^c}(T)=N(T)\setminus Q$.

Consider a fixed pair of subsets $Q_1,Q_2\subset [k]$. Fixing, for
the moment, the random variables $\bX_Q$, let us consider the
expectation of, say, $B(P_{Q^c}(T_1),P_{Q^c}(T_2))$. For $T\subset
[k]$, let ${\bX}_T$ be the collection of random variables $X_i$
with $i\in T$, and let
\[
{\mathcal P}(\bX_T)=\ \Bigl\{x\in [0,1]\colon \sum_{i\in T}
U(x,X_i)>0\Bigr\}.
\]
Then $P_{Q^c}(T_i)=\{j\in {Q^c}\colon X_j\in {\mathcal P}(\bX_{T_i})\}$, and
\[
\begin{aligned}
\E_{\bX_{Q^c}}\Bigl[B(P_{Q^c}(T_1),P_{Q^c}(T_2))\Bigr] &=
|Q^c|(|Q^c|-1)\int_{{\mathcal P}(\bX_{T_1})\times {\mathcal
P}(\bX_{T_2})} U(x,y)\,dx\,dy
\\
&\quad + |Q^c|\int_{{\mathcal P}(\bX_{T_1})\cap{\mathcal
P}(\bX_{T_2})} U(x,x)\,dx
\\
&\leq k^2\|U\|_\square +k.
\end{aligned}
\]
It is not hard to see that the random variable
$Y=B(P_{Q^c}(T_1),P_{Q^c}(T_2))$ is highly concentrated.  Indeed,
consider $Y$ as a function of the random variables in $\bX_{Q^c}$. If
we change one of these variables, then $Y$ changes by at most $4k$,
implying that with probability at least {$1-e^{-k\rho^2/32}$, $Y\leq
\E(Y)+\rho k^2$.} The same bound holds for the random variable
$\tilde Y=-B(N_{Q^c}(T_1),N_{Q^c}(T_2))$. As a consequence, the
maximum in \eqref{B+bound} obeys this bound with probability at least
$1-2^{2p+1}{ e^{-k\rho^2/32}}$. Since $Y,\tilde Y\leq k^2$ for all
$\bX$, we conclude that
\[
\E\Bigl[\|\bH(k,U)\|_\square \Bigr]= \frac
1{k^2}\E\Bigl[\|B\|_\square\Bigr] \leq \|U\|_\square+\frac {2}{\sqrt
p} + \frac{4p}k +\frac 1k + {\rho +2^{2p+1}e^{-k\rho^2/32}.}
\]
Choosing $p$ and $\rho$ of the form $p=\lceil\alpha \sqrt k\rceil$
and $\rho=\beta k^{-1/4}$ with $\alpha=(4\sqrt{\log 2})^{-1}$ and
$\beta=4(\log 2)^{1/4}+4/10$, this implies that for $k\geq 10^4$, we
have
\begin{equation}
\begin{aligned}
\E\Bigl[\|\bH(k,U)\|_\square \Bigr] &\leq\|U\|_\square+ \frac
{1}{k^{1/4}}\Bigl(8(\log 2)^{1/4}+0.534\dots)\Bigr) \leq
\|U\|_\square+\frac {8}{k^{1/4}},
\end{aligned}
\end{equation}
which is the upper bound in \eqref{ExpH-UpBd}.

To prove concentration, we use that $\|\bH(k,U)\|_\square$ changes by
at most $\frac 4k\|U\|_\infty$
if we change one of the random variables $X_1,\dots,
X_k$, so applying Lemma \ref{lem:Azuma} (i)
we get that its values are
highly concentrated around its expectation:
\[
\Prob\Bigl(\Bigl|\|\bH(k,U)\|_\square -
\E\Bigl[\|\bH(k,U)\|_\square\Bigr]\Bigr|
>{\frac{2}{k^{1/4}}}\Bigr)
< 2e^{-\sqrt k/8}.
\]
This completes the proof of (i).

(ii)
We couple
$\bG(k,U_1)$ and $\bG(k,U_2)$ as follows: as in the proof
of (i)
we chose $X_1,\dots,X_k$ to be i.i.d., distributed uniformly
in $[0,1]$.  In addition, we independently choose $k(k+1)/2$ random
variables $Y_{ij}=Y_{ji}$
uniformly
at random from $[0,1]$.  In terms of these random variables,
we then define $G_1$ to be the graph on $[k]$
which has an edge between $i$ and $j$
whenever $U_1(X_i,X_j)<Y_{ij}$, and
$G_2$ to be the graph which has an edge
between $i$ and $j$
whenever $U_2(X_i,X_j)<Y_{ij}$.
Then $G_1$, $G_2$ is a coupling of
$\bG(k,U_1)$ and $\bG(k,U_2)$, and
\[
d_\square(G_1,G_2)=
\frac{1}{k^2}\max_{S_1,S_2\subset[k]}|B(S_1,S_2)|,
\]
where $B$ is the matrix with entries
$
B_{ij}=\one_{Y_{ij}<U_1(X_i,X_j)}
-\one_{Y_{ij}<U_2(X_i,X_j)}
$.
We again have to
bound the expectation of
$B(P_{Q^c}(T_1),P_{Q^c}(T_2))$.
As before, we fix the variables
$X_i$ for $i\in Q$, but now we also fix all random variables
$Y_{ij}$ for which $\{i,j\}$ intersects $Q$. In order to
calculate expectations with
respect to the remaining random variables,
we express $B(P_{Q^c}(T_1),P_{Q^c}(T_2))$
as the sum
$\sum_{i,j\in Q^c} B_{i,j}^+(T_1,T_2)$, where
\[
B_{i,j}^+(T_1,T_2)
=B_{ij}\one_{i\in P(T_1)}\one_{j\in P(T_1)}.
\]
For $i\in Q^c$, the event that
$i\in P(T)$ for some $T\subset Q$
can then be re-expressed as the event
that $X_i$ lies in the set
\[
\Bigl\{x\in [0,1]\colon \sum_{j'\in T}
\one_{Y_{ij'}<U_1(x,X_{j'})}
>\sum_{j'\in T}\one_{Y_{ij'}<U_2(x,X_{j'})}\Bigr\}.
\]
Observing that this set only depends on the fixed random
variables, we can now proceed as before to calculate the expectation
of $B_{i,j}^+(T_1,T_2)$, and hence the expectation
of $B(P_{Q^c}(T_1),P_{Q^c}(T_2))$.  This, together with a similar
(again much easier)
calculation for the lower bound, leads to the
estimate
\begin{equation}
\label{ExpG-UpBd}
-\frac 2k
\leq \E\Bigl[d(G_1,G_2)_\square \Bigr]-\|U_1-U_2\|_\square
\leq
\frac{8}{k^{1/4}},
\end{equation}
as before valid for $k\geq 10^4$. Concentration is again proved with
the help of Lemma~\ref{lem:Azuma}.
\end{proof}

Essentially the same proof also gives Theorem~\ref{thm:dist-test}:

\medskip

\noindent{\bf Proof of Theorem~\ref{thm:dist-test}} We generate the
set $S$ by choosing $X_1,\dots,X_k$ uniformly without replacement
from $V$.  In this way, we rewrite $d_\square(G_1[S],G_2[S])$ in
terms of the matrix $B=B(\mathbf X)$ with entries
$B_{ij}=\beta_{X_iX_j}(G_1)-\beta_{X_iX_j}(G_2)$. Observing that
\[
|B_{ij}|\leq C=2\max\{\|G_1\|_\infty,{\|G_2\|_\infty}\},
\]
we may proceed exactly as in the
proof of Theorem~\ref{thm:NormSample} (i),
leading to the bound
\begin{equation}
\label{ExpG[S]-UpBd}
-\frac {2C}k
\leq \E[d_\square(G_1[S],G_2[S]]-d_\square(G_1,G_2)
\leq \frac {8C}{k^{1/4}}.
\end{equation}
Using Lemma~\ref{lem:Azuma} (ii) to prove concentration, we get
the bound \eqref{ProbG[X]-G}.
\hfill$\square$

We now come to the main theorem about sampling.

\begin{theorem}\label{thm:sample2}
Let $k$ be a positive integer.

\smallskip

{\rm (i)} If  $U\in\WW$, then with probability at least
$1-e^{-k^2/(2\log_2 k)}$, we have
\[
\delta_\square(U, \bH(k,U))
\leq \frac{10}{\sqrt{\log_2 k}}\|U\|_\infty.
\]

\smallskip

{\rm (ii)} If  $U\in\WWz$, then with probability at least
$1-e^{-k^2/(2\log_2 k)}$, we have
\[
\delta_\square(U, \bG(k,U))
\leq \frac{10}{\sqrt{\log_2 k}}.
\]
\end{theorem}

\begin{proof}
We again first bound expectations and then prove concentration, and
as before, we assume without loss of generality that
$\|U\|_\infty\leq 1$. Finally, we may assume that $k\geq 2^{25}\geq
10^4$, since otherwise the bounds of the theorem are trivial.

{In a first step, we use the weak Szemer\'edi approximation in Lemma
\ref{lem:W-szem} and the sampling bound \eqref{ExpH-UpBd} to show
that it is enough to consider the case where $U$ is a step function.
Indeed, given $\eps>0$, let $U_1$ be a step function  with $q\leq
4^{1/\eps^{2}}$ steps such that
\begin{equation}\label{U-U1}
\|U-U_1\|_\square \leq \eps.
\end{equation}
By the bound \eqref{ExpH-UpBd},
we have that
\[
\E\Bigl[\delta_\square(\bH(k,U_1), \bH(k, U))\Bigr]
\leq \E[\,\|\bH(k,U-U_1)\|_\square]\leq\eps
+\frac{16}{k^{1/4}}.
\]
As a consequence,
\begin{equation}
\label{step-reduction}
\begin{aligned}
\E\Bigl[\delta_\square(U, \bH(k,U))\Bigr]
&\leq \delta_\square(U, U_1) +
\E\Bigl[\delta_\square(U_1, \bH(k,U_1))\Bigr]
+\E\Bigl[ \delta_\square(\bH(k,U_1), \bH(k,U))\Bigr]
\\
&\leq
2\eps+\frac{16}{k^{1/4}}+\E\Bigl[\delta_\square(U_1, \bH(k,U_1))\Bigr].
\end{aligned}
\end{equation}

We are thus left with the problem of sampling from the step function $U_1$.}
Let $U_1$ have steps $J_1,\dots,J_q\subseteq[0,1]$, and
$\lambda(J_i)=\alpha_i$. { Let $X_1,\dots,X_k$ be independent random
variables that are uniformly distributed on $[0,1]$, and let $Z_i$ be
the number of points $X_j$ that fall into the set $J_i$.}
It is easy to compute that
\[
\E(Z_i)=\alpha_i k, \qquad \Var(Z_i) = (\alpha_i-\alpha_i^2)k <
\alpha_ik.
\]
Construct a partition of $[0,1]$ into measurable sets
$J_1',\dots,J'_q$ such that { $\lambda(J_i')=Z_i/k$ and}
\[
\lambda(J_i\cap J'_i)=\min(\alpha_i, Z_i/k),
\]
and also construct a symmetric function $U'\in\WW$ such that the
value of $U'$ on $J_i'\times J_j'$ is the same as the value of $U_1$
on $J_i\times J_j$. Then $U'$ is a step function representation of
$\bH(k,U_1)$, and it agrees with $U_1$ on the set
$
Q=\cup_{i,j=1}^q (J_i\cap J'_i)\times (J_j\cap J'_j).
$
Thus
\begin{align*}
\delta_\square(U_1,\bH(k,U_1)) &\le \|U_1-U'\|_\square \le \|U_1-U'\|_1
\le
{ 2\big(} 1-\lambda(Q) { \big)} = { 2\Bigl(} 1-\Bigl(\sum_i
\min\Bigr(\alpha_i, \frac{Z_i}{k}\Bigl)\Bigr)^2 { \Bigr)}
\\
&\le { 4}\Bigl(1-\sum_i \min\Bigr(\alpha_i, \frac{Z_i}{k}\Bigl)\Bigr)
= { 2}\sum_i\Bigl|\alpha_i-\frac{Z_i}{k}\Bigr| \le { 2\bigg(}q
\sum_i\Bigl(\alpha_i-\frac{Z_i}{k}\Bigr)^2{\bigg)^{1/2}},
\end{align*}
which we rewrite as
\[
\Bigl(\delta_\square(U_1,\bH(k,U_1)\Bigr)^2\leq 4
q\sum_i\Bigl(\alpha_i-\frac{Z_i}{k}\Bigr)^2.
\]
The expectation of the right hand side is
\[
\frac{{ 4}q}{k^2}\sum_i \Var(Z_i)<\frac{{ 4}q}{k},
\]
so by Cauchy--Schwarz
\begin{equation}
\label{STEP-SAMPLE}
\E\Bigl[\delta_\square(U_1,\bH(k,U_1))\Bigr]
\leq\sqrt{\frac{{ 4}q}{k}}.
\end{equation}
Inserted into \eqref{step-reduction}, this gives
\[
\E\Bigl[\delta_\square(U, \bH(k,U))\Bigr]
\leq
2\eps+\frac{16}{k^{1/4}}+
\sqrt{\frac{{ 4}n}{k}}
\leq 2\eps+\frac{16}{k^{1/4}}+
\frac 2{k^{1/2}}2^{1/\eps^{2}}.
\]
Choosing $\eps=2/(\log_2 k)$ and recalling that $k\geq 2^{25}$, this
gives the upper bound
\begin{equation}
\label{E-sample-bound}
\E\Bigl[\delta_\square(U, \bH(k,U))\Bigr]
\leq \frac{1}{\sqrt{\log_2 k}}
\Bigl(4+ (16+2)
\frac{\sqrt{\log_2 k}}{k^{1/4}}\Bigr)
\leq \frac{6}{ \sqrt{\log_2 k}}.
\end{equation}
Observing that $\delta_\square(U, \bH(k,U))$ changes by at most
$4/k$ if one of the random variables $X_i$ changes its value,
we finally use Lemma~\ref{lem:Azuma} (i)
to complete the proof of the first statement.

(ii) The proof of the
second statement is completely analogous: we first show that
\begin{equation}
\label{E-G-sample-bound}
\E\Bigl[\delta_\square(U, \bG(k,U))\Bigr]\leq
\frac{1}{\sqrt{\log_2 k}}
\Bigl(4+ (8+2)
\frac{\sqrt{\log_2 k}}{k^{1/4}}\Bigr)
\leq\frac{5}{\sqrt{\log_2 k}}
,
\end{equation}
and then use Lemma~\ref{lem:Azuma} (i) to prove concentration.
\end{proof}

The proof again generalizes
to samples from weighted graphs, this time leading to
Theorem~\ref{thm:sample}.

\noindent{\bf Proof of Theorem~\ref{thm:sample}.}
{ Again, the proof is analogous to the proof of
statement (i) above, except that we now use the original Frieze-Kannan Lemma
(Lemma~\ref{lem:szem-weak}) for graphs instead of our Weak
Szemer\'edi Lemma (Lemma~\ref{lem:W-szem}) for graphons.  One also
needs to generalize the bound \eqref{STEP-SAMPLE} to samples
$X_1,\dots,X_k$ chosen uniformly without replacement from $V(G)$, but
this is again straightforward: the random variable $Z_i$ is now the
number of points that fall into the class $V_i$ of a weak Szemer\'edi
partition $\PP=(V_1,\dots,V_q)$, and its variance is bounded by
$\alpha_i k$, where $\alpha_i=|V_i|/|V(G)|$. Continuing as in the
above proof, these considerations now lead to the bound
\begin{equation}
\label{E-G[X]-sample-bound}
\E\Bigl[\delta_\square(G,\RandG(k,G))\Bigr]\leq
\frac{1}{\sqrt{\log_2 k}} \Bigl(4+ (16+2) \frac{\sqrt{\log_2
k}}{k^{1/4}}\Bigr) \leq\frac{6}{\sqrt{\log_2 k}},
\end{equation}
where we again assume that the weights have been rescaled so that
$\|G\|_\infty=1$. Using Lemma~\ref{lem:Azuma} (ii) to show
concentration, this gives the bound \eqref{sample-bd}.
}
\hfill$\square$

\subsection{Proof of Theorem~\ref{W-CONV-CHAR}}
\label{sec:proof-main-W}

\subsubsection{Proof of Theorem~\ref{W-CONV-CHAR} (a)}
\label{HOMDIST}

Let $V(F)=[n]$ and $E(F)=\{e_1,\dots,e_m\}$. For $t=1,\dots, m$, let
$i_t$, $j_t$ be the endpoints of the edge $e_t$, and let
$E_t=\{e_1,\dots,e_t\}\subset E(F)$. Then $t(F,W)-t(F,W')$ can be
rewritten as
\begin{align}
t(F,&W)-t(F,W')=\int\limits_{[0,1]^n} \Bigl(\prod_{ij\in E(F)}
W(x_i,x_j)-\prod_{ij\in E(F)} W'(x_i,x_j)\Bigr)\,dx_1\dots dx_n
\nonumber\\
&=\sum_{t=1}^{m} \int\limits_{[0,1]^n} \prod_{s<t} W(x_{i_s},x_{j_s})
\prod_{s>t}W'(x_{i_s},x_{j_s}) \Bigl(W(x_{i_t},x_{j_t})-
W'(x_{i_t},x_{j_t})\Bigr)\,dx_1\dots dx_n. \label{TELESCOPE}
\end{align}

Take any term in this sum, and for notational convenience, assume
that $i_t=1$ and $j_t=2$. Let $X(x_1,x_3,\dots,x_n)$ be the product
of those factors in $\prod_{s<t} W(x_{i_s},x_{j_s}) \prod_{s>t}
W'(x_{i_s},x_{j_s})$ that contain $x_1$, and let $Y(x_2,\dots,x_n)$
be the product of the rest. Then we have
\begin{align*}
\int\limits_{[0,1]^n}& \prod_{s<t} W(x_{i_s},x_{j_s})
\prod_{s>t}W'(x_{i_s},x_{j_s}) \Bigl(W(x_{i_t},x_{j_t})-
W'(x_{i_t},x_{j_t})\Bigr)\,dx_1\dots dx_n \\ &=
\int\limits_{[0,1]^{n-2}} \Bigl(\int\limits_{[0,1]^2}
X(x_1,x_3\dots,x_n)Y(x_2,\dots,x_n) \bigl(W(x_1,x_2)-
W'(x_1,x_2)\bigr)\,dx_1\,dx_2\Bigr) dx_3 \dots dx_n.
\end{align*}
Here the interior integral is bounded by
\begin{align*}
\Bigl|\int\limits_{[0,1]^2}  &X(x_1,x_3\dots,x_n)Y(x_2,\dots,x_n)
\bigl(W(x_1,x_2)-
W'(x_1,x_2)\bigr)\,dx_1\,dx_2\Bigr|
\\&
\le \|X\|_\infty \|Y\|_\infty\|W-W'\|_{\infty\to 1}.
\end{align*}
Substituting into \eqref{TELESCOPE}, and
using that
\[
\|X\|_\infty\|Y\|_\infty \leq
\|W\|_\infty^{t-1}\|W'\|_\infty^{m-t}\le C^{m-1}
\]
and (by \eqref{CUT-IO-W})
\[
\|W-W'\|_{\infty\to 1}\le 4 \|W-W'\|_\square,
\]
we get
\[
|t(F,W)-t(F,W')|\le 4 mC^{m-1}\|W-W'\|_\square,
\]
Using the representation in Lemma \ref{lem:delta=delta} for the
$\delta_\square$ distance and the fact that $t(F,W)=t(F,W^\phi)$
whenever $\phi$ is a measure-preserving function from $[0,1]$ to
$[0,1]$, this bound implies the bound \eqref{LEFTDIST-delta}.
$\hfill\square$

\begin{remark}
The above proof can easily be generalized to show
that \begin{equation}
\label{LEFTNORM} |t(F,W)|\le { 4 \|W\|_\infty^{|E(F)|-1}}
\|W\|_\square.
\end{equation}
for all $W\in\WW$ and all simple graphs $F$.
Also, it is not hard to show
that the factor $4$ in \eqref{LEFTDIST-delta} and \eqref{LEFTNORM}
is not needed if $W$ and $W'$ are non-negative.
\end{remark}

\subsubsection{Proof of Theorem~\ref{W-CONV-CHAR} (b)}
\label{sec:left-vs-metric-conv}

Without loss of generality, we may assume that $\|W\|_\infty,
\|W'\|_\infty\leq 1$ (otherwise we just rescale both $W$ and $W'$ by
the maximum of these two numbers). Let $U,U'\in\WWz$ be the graphons
$U=\frac 1{2}W+\frac 12$ and $U'=\frac 1{2}W'+\frac 12$. Then
$\delta_\square(W, W')=2\delta_\square(U,U')$, so it is enough to
prove that $\delta_\square(U,U')\leq{ 11/\sqrt{\log_2 k}}$. We will
prove this bound by relating the distance of $U$ and $U'$ to the
distance of the random graphs $\bG(k,U)$ and $\bG(k,U')$.  To this
end, we need an expression for the probability that $\bG(k,U)$ is
equal to some given graph $F$ on $k$ nodes. We first use the
relations \eqref{W-random-Exp} to express $t(F,U)$ as a sum over
graphs $G$ on $k$ nodes.  Combined with \eqref{T0T1} and the fact
that for all graphs $F'$ and $G$ on $k$ nodes, $t_\ind(F',G)=0$
unless $G$ is isomorphic to $F'$, we have that
\[
\begin{aligned}
t(F,U)&=\sum_{G}\Prob\Bigl(\bG(k,U))=G\Bigr)t_\inj(F,G)
\\
&=\sum_{G}\sum_{F'\supset F} \Prob\Bigl(\bG(k,U))=G\Bigr)t_\ind(F',G)
\\
&=\frac 1{k!}\sum_{F'\supset F} \Prob\Bigl(\bG(k,U))=F'\Bigr)
\Bigl(\ind(F',F')\Bigr)^2,
\end{aligned}
\]
{where $\ind(\cdot,\cdot)$ is defined in \eqref{ind-def}.}
With the help of inclusion-exclusion, this leads to
\[
\begin{aligned}
\Prob\Bigl(\bG(k,U))=F\Bigr) =k!\sum_{F'\supset F}
(-1)^{|E(F')\setminus E(F)|} \Bigl(\ind(F',F')\Bigr)^{-2}t(F',U),
\end{aligned}
\]
which in turn implies that
\[
\begin{aligned}
&\sum_F\Bigl|\Prob\Bigl(\bG(k,U))=F\Bigr)-\Prob\Bigl(\bG(k,U'))=F\Bigr)\Bigr|
\leq k!\sum_{\bsa F,F':\\ E(F')\supset E(F)\esa}
\Bigl|t(F',U)-t(F',U')\Bigr|
\end{aligned}
\]
where the sum runs over graphs $F$ and $F'$ on $k$ nodes. To
continue, we need to relate the homomorphism densities of $U$ and
$U'$ to those of $W$ and $W'$. To this end, we insert the relation
$U=\frac 12\bigl(W+1\bigr)$ into the definition of $t(F',U)$. For a
graph $F'$ on $k$ nodes, this leads to the identity
\[
t(F',U)=2^{-|E(F')|}\sum_{F''\subset F'} t(F'',W)
\]
where the sum goes over all subgraphs  that have the same node set as
$F'$.  Using the assumption of the theorem, we thus obtain the bound
$ \bigl|t(F',U)-t(F',U')\bigr| \leq 3^{-k^2} $ which in turn implies
that
\[
\sum_F\Bigl|\Prob\Bigl(\bG(k,U))=F\Bigr)-\Prob\Bigl(\bG(k,U'))=F\Bigr)\Bigr|
\leq k!\sum_{\bsa F,F':\\ E(F)\subset E(F')\esa} 3^{-k^2}
=  k!3^{-k^2}3^{k(k-1)/2}.
\]
Bounding
$k!$ (rather crudely) by $3^{k^2/2}$, we note that the right
hand side is smaller than $3^{-k/2}$.

As a consequence, $\bG(k,U)$ and $\bG(k,U')$ can be coupled in such a
way that $\bG(k,U)=\bG(k,U')$ with probability at least $1-3^{-k/2}$.
Combined with the triangle inequality and
the bound \eqref{E-G-sample-bound}, we obtain
that
\[
\begin{aligned}
\delta_\square(U,U') &\leq
\E{\delta_\square(\bG(k,U),\bG(k,U'))}
+\E[\delta_\square(\bG(k,U),U)]+
\E[\delta_\square(\bG(k,U'),U')]
\\
&\leq 3^{-k/2}+{\frac{10}{\sqrt{\log_2 k}}}
\leq
{\frac{11}{\sqrt{\log_2 k}}.}
\end{aligned}
\]
\hfill$\square$

\section{Convergence in Norm and Uniform Szemer\'edi Partitions}
\label{sec:norm+szem}

\subsection{Comparison of Fractional and Non-Fractional Overlays}
\label{sec:compdist}

Let $G$ and $G'$ be two graphs on $n$ nodes, both with nodeweights
$1/n$.   Consider any labeling that attains the minimum in the
definition \eqref{delta-hat} of $\iv\delta_\square$, and identify
the nodes of $G$ and $G'$ with the same label. In this case, we say
that $G$ and $G'$ are {\it optimally overlaid.}

In addition to this distance we also defined the distance
$\delta_\square$, given in terms of fractional overlays $X\in
\XX(G,G')$, see \eqref{delta-altdef}. Since every bijection between
the nodes of $G_1$ and $G_2$ defines a fractional overlay $X$, we
trivially have that
\[
\delta_\square(G,G')\le \iv\delta_\square(G,G').
\]
This inequality can be strict. Let $G=K_2$, and let $G'$ be a graph
with two nonadjacent nodes but with a loop at each node. It is easy
to see that $\iv\delta_\square(G,G')=1/4$, but
$\delta_\square(G,G')=1/8$ (the best fractional overlay is
$X_{iu}=1/4$ for all $i\in V(G)$ and $u\in V(G')$). To give an
example without loops, let $G=K_{3,3}$, and let $G'$ consist of two
disjoint triangles $\Delta_1$ and $\Delta_2$. There are only two
essentially different ways to overlay these graphs; the better one
maps two nodes of $\Delta_i$ into one color class of $K_{3,3}$ and
the third one into the other color class. The number of edges in
$G'$ between $\Delta_1$ and $\Delta_2$ is $5$, whence
$\iv\delta_\square(G,G')\ge 5/36$. (One can check that equality
holds.) On the other hand, let us double each node in both graphs,
to get $G(2)=K_{6,6}$ and $G'(2)=\Delta_1(2)\cup \Delta_2(2)$. Let
us map one copy of each twin node of $G'(2)$ into one color class of
$K_{6,6}$, and its pair into the other color class. Case distinction
shows that the worst choice for $S$ and $T$ in the definition of
$d_\square$ is $S=V(\Delta_1(2))$ and $T=V(\Delta_2(2))$, showing
that
\[
\delta_\square(G,G')\le \iv\delta_\square(G(2),G'(2))=\frac{18}{144}
= \frac{1}{8}<\frac{5}{36}.
\]

We have no example disproving the possibility that
$\iv\delta_\square(G_1,G_2)=O(\delta_\square(G_1,G_2))$, but we are
only able to prove a much weaker inequality given in Theorem
\ref{TWIN-SQUARE}. We start with a simple but very weak bound we
will need.

\begin{lemma}\label{EASY-DELTA}
Let $G_1$ and $G_2$ be weighted graphs on
$n$ nodes.  If both $G_1$ and $G_2$ have
nodeweights one, then
\[
\iv\delta_\square(G_1,G_2)\le n^{{ 6}} \delta_\square(G_1,G_2).
\]
\end{lemma}

\begin{proof}
Let $(X_{ui})$ be an optimal fractional overlay of $G_1$ and $G_2$,
{normalized in such a way that
$\sum_i X_{ui}=\sum_v X_{vj}=1/n$ for all $n$.}
We claim that there is a bijection $\pi:~V(G_1)\to V(G_2)$ such that
$X_{u\pi(u)}\ge 1/n^{{3}}$ for all $u\in V(G_1)$. This follows
from a routine application of the Marriage Theorem: if there is no
such bijection, then there are two sets $S\subseteq V(G_1)$ and
$T\subseteq V(G_2)$ such that $|S|+|T|>n$ and $X_{st}<1/n^{{3}}$
for all $s\in S$
{and all $t\in T$.} But then
\begin{align*}
\frac{|S|}{n}&=\sum_{u\in S}\sum_{i\in V(G_2)} X_{iu} =\sum_{u\in
S}\sum_{i\in T} X_{iu} +\sum_{u\in S}\sum_{i\in V(G_2)\setminus T}
X_{iu}\\ &\le \frac{1}{n^3}|S|\cdot|T| +\frac{|V(G_2)\setminus T|}{n}
< \frac{1}{n}+\frac{|S|-1}{n},
\end{align*}
a contradiction. Thus a map $\pi$ with the desired properties exists.

Let $G_1'$ be the image of $G_1$ under this map.
Then
\begin{align*}
\iv\delta_\square(G_1,G_2)&
\leq\max_{s\in V(G_1)}\max_{t\in  V(G_2)}
|\beta_{st}(G_1)-\beta_{\pi(s)\pi(t)}(G_2)|
\\
&\le \max_{s\in V(G_1)}\max_{t\in  V(G_2)}
n^6X_{s\pi(s)}X_{t\pi(t)}|\beta_{st}(G_1)
-\beta_{\pi(s)\pi(t)}(G_2)|\\
&\le n^6\delta_\square(G_1,G_2),
\end{align*}
which proves the lemma.
\end{proof}

Now we are able to prove Theorem~\ref{TWIN-SQUARE}, which shows that
the two distances $\delta$ and $\iv\delta$ define the same topology.
As the reader may easily verify, the proof below gives an exponent
of $1/3$ instead of the exponent of $1/{67}$ from
Theorem~\ref{TWIN-SQUARE} if the number of nodes is large enough.
However, even under this assumption, we could not obtain a linear
bound, i.e., bi-Lipschitz equivalence.

\medskip

{\bf\noindent Proof of Theorem~\ref{TWIN-SQUARE}} The first
inequality, as remarked before, is trivial. To prove the second,
write $\delta_\square(G_1,G_2)^{1/{ 67}}=\eps$. If $n\leq
\eps^{-11}$, then the bound follows by Lemma \ref{EASY-DELTA}, and
if $\eps\geq {2/36}$, the bound is trivial, so we may assume that
\begin{equation}\label{N-BOUND}
n>\eps^{-11}{\qquad\text{and}\qquad \eps\leq \frac 1{16}.}
\end{equation}
Consider an optimal overlay of $W_{G_1}$ and $W_{G_2}$ (in other
words, consider a measure-preserving bijection
$\phi:~[0,1]\to[0,1]$) such that
\[
\|W_{G_1}-W_{G_2}^\phi\|_\square = \delta_\square(G_1,G_2)
=\eps^{{67}}.
\]
Let us select a set $Z$ of $k=\lceil n/\eps\rceil$ random points from
$[0,1]$. Let $H_1=W_{G_1}[Z]$ and $H_2=W_{G_2^\phi}[Z]$. Then by
Lemma \ref{thm:NormSample}, we get that with probability
at least
{ $1-e^{-\frac 18\sqrt{n/\eps}}\geq 1-8\eps^6$,}
\[
d_\square(H_1,H_2)\le \|W_{G_1}-W_{G_2}^\phi\|_\square +
\frac{{20}}{k^{1/4}} = \eps^{ 67} +\frac{{20}}{k^{1/4}}.
\]
Each element $z\in Z$ corresponds to a node $i_z\in V(G_1)$ and a
node $j_z\in V(G_2)$. These pairs $(i_z,j_z)$ form a bipartite graph
$F$ with color classes $V(G_1)$ and $V(G_2)$ (which we assume are
disjoint).

\begin{claim}\label{MATCHING}
With probability at least $1-{(2/e)^{-2n}}$, the bipartite graph $F$ has a
matching of size at least $(1-2\eps)n$.
\end{claim}

To prove this claim, we use K\"onig's Theorem: if $F$ does not
contain a matching of size $(1-2\eps)n$, then its edges can be
covered by a set $X$ of nodes with $|X|<(1-2\eps)n$. Let
$Y_i=V(G_i)\setminus X$.  Then there is no edge of $F$ between $Y_1$
and $Y_2$.

On the other hand, $|Y_1|+|Y_2|\ge (1+2\eps)n$. Let
$J_i\subseteq[0,1]$ be the union of intervals representing $Y_i$ in
$W_{G_i}$, so that $\lambda(J_1)+\lambda(J_2)\ge 1+2\eps$, and hence
also $\lambda(J_1)+\lambda(\phi(J_2))\ge 1+2\eps$. This implies that
$\lambda(J_1\cap \phi(J_2))\ge2\eps$. The random set $Z$ avoided this
intersection; the probability of this happening is at most
\[
(1-2\eps)^k < e^{-2\eps k}\leq e^{-2n}.
\]
Since there are at most $4^n$ pairs of sets $(Y_1,Y_2)$, the
probability that $F$ does not have a matching of cardinality at least
$(1-2\eps)n$ is less than $4^n e^{-2n}$. This proves
the claim.

Now let $(i_1,j_1)\dots, (i_m,j_m)$ be a maximum matching in $F$.
With positive probability, we have both $d_\square(H_1,H_2)\le
\eps^{ 67}+{ 20}/k^{1/4}$ and
$m\ge (1-\frac{2n}{k})n$. Fix a choice of
$Z$ for which this happens. Let $(i_{m+1},j_{m+1}),\dots,(i_n,j_n)$
be an arbitrary pairing of the remaining nodes of $V(G_1)$ and
$V(G_2)$. We claim that the pairing $\pi:~i_r\mapsto j_r$ gives an
overlay of $G_1$ and $G_2$ with small $d_\square$ distance.

Let $S,T\subseteq V(G_2)$, and let $S'=S\cap\{j_1,\dots,j_m\}$,
$T'=T\cap\{j_1,\dots,j_m\}$. Then $|S\setminus S'|\le\frac{2n^2}{k}$,
and $|T\setminus T'|\le\frac{2n^2}{k}$, and hence
\[
|e_{G_1'}(S,T)-e_{G_2}(S,T)| \le |e_{G_1'}(S',T')-e_{G_2}(S',T')|
+2\Bigl(2 n \frac{2n^2}{k}+\Bigl(\frac{2n^2}{k}\Bigr)^2\Bigr).
\]
Here
\begin{align*}
|e_{G_1'}(S',T')-e_{G_2}(S',T')| &=  |e_{H_1}(S',T')-e_{H_2}(S',T')|
\le d_\square(H_1,H_2)k^2
\le (\eps^{ 67}+\frac{{20}}{k^{1/4}}) k^2,
\end{align*}
and hence
\[
\begin{aligned}
\frac{|e_{G_1'}(S,T)-e{G_2}(S,T)|}{n^2}
&\le
{\frac{8n}k+\frac{8n^2}{k^2}}
+\frac{\eps^{ 67} k^2}{n^2} + \frac{{ 20}k^{\frac{7}{4}}}{n^2}
\\
&\le
{8\eps+8\eps^2 +
    \eps^{67}(1+\eps^{-1})^2 + 20(1+\eps^{-1})^{7/4}
    n^{-1/4}
    }
\\
&{\leq 8\eps+8\eps^2 +
    \eps^{65}(1+\eps)^2 + 20\eps(1+\eps)^{7/4}
    \leq 32\eps
    }
\end{aligned}
\]
where we used \eqref{N-BOUND} in the last two bounds.
\hfill$\square$

\subsection{Convergence in Norm}
\label{sec:Norm-Conv}

Let $(G_n)$ be a convergent sequence of weighted graphs.
Theorem~\ref{LEFTCONV=DELTACONV} then implies that there  exists a
graphon $W\in \WW$ such that $W_{G_n}\to W$ in the $\delta_\square$
distance. This does not imply, however, that $W_{G_n}\to W$ in the
$\|\cdot\|_\square$ norm. It turns out, however, that the graphs in
the sequence $(G_n)$ can be relabeled in such a way that this becomes
true, provided $(G_n)$ has {\it no dominant nodes} in the sense
that
\[
\frac {\max_i\alpha_i(G_n)}{\alpha_{G_n}}\to 0
\]
as $n\to\infty$. This is the content of the following lemma, which
will be useful when discussing testability.

\begin{lemma}
\label{lem:norm-conv-graph} Let $(G_n)$ be a sequence of weighted
graphs with uniformly bounded edgeweights, {and no dominant nodes}.
If
\[
\delta_\square(U,W_{G_n})\to 0
\]
for some $U\in\WW$, then the graphs in the sequence $(G_n)$  can be
relabeled in such a way that the resulting sequence $(G_n')$ of
labeled graphs converges to $U$ in the cut-norm:
\[
\|U-W_{G_n'}\|_\square\to 0.
\]
\end{lemma}

\begin{proof}
{We first prove the lemma for graphs with nodeweights one.} Let
$m(n)=|V(G_n)|$, and let $\PP_{m(n)}$ be a partition of $[0,1]$ into
consecutive intervals of length $1/m(n)$. By
Lemma~\ref{lem:as-approx}, we have that
$\|U-U_{\PP_{m(n)}}\|_\square\to 0$, so combined with the assumption
that $\delta_\square(U,W_{G_n})\to 0$ we conclude that
$\delta_\square(U_{\PP_{m(n)}},W_{G_n})\to 0$.  But the left hand
side can be expressed as the distance of two weighted graphs on the
same number of nodes,
$\delta_\square(U_{\PP_{m(n)}},W_{G_n})=\delta_\square(U/\PP_{m(n)},G_n)$,
so by Theorem~\ref{TWIN-SQUARE}, we get that
\[
\hat \delta_\square(U/\PP_{m(n)}, G_n)\to 0.
\]
This means that the graphs in the sequence $(G_n)$ can be relabeled
in such a way that for the resulting graph sequence, $(G_n')$, we have
\[
\|U_{\PP_{m(n)}}-W_{G_n'}\|_\square
=d_\square(U/{\PP_{m(n)}},{G_n'})\to 0.
\]
Combined with the fact that
$\|U-U_{\PP_{m(n)}}\|_\square\to 0$,
this gives the statement of the lemma
{for graphs with nodeweights one.

To prove the lemma for general sequences without dominant nodes, we
use the Weak Regularity Lemma to approximate $(G_n)$ by a sequence
of graphs $(\tilde G_n)$ with nodeweights one. Indeed, let us assume
without loss of generality that all graphs in the sequence $(G_n)$
have total nodeweight $\alpha_{G_n}=1$.  Define
$a_n=\max_i\alpha_i(G_n)$, and choose $\eps_n$ in such a way that
$\eps_n\to 0$ and $a_n 2^{40/\eps^2_n}\to 0$ as $n\to\infty$. With
the help of Corollary~\ref{cor:unif-Sem} (ii), we then construct a
partition $\PP_n$ of $V(G_n)$ into $q_n\leq 2^{20/\eps_n^2}$ classes
such that $d_\square(G_n, (G_n)_{\PP_n})\to 0$ and the classes in
$\PP_n$ have almost equal weights (in the sense that $|\sum_{x\in
V_i}\alpha_x(G_n)-\sum_{y\in V_j}\alpha_y(G_n)|\leq a_n$ for all
$i,j\in[q_n]$). {Consider the sequence of graphs $\tilde G_n$ that
are obtained from $G_n/\PP_n$ by changing all nodeweights to $1$.
Since} the classes of $\PP_n$ have almost equal weights, we have
that $\|W_{\tilde G_n}-W_{G_n/\PP_n}\|_\square\leq q_n^2 a_n\to 0$
which in turn implies that $\delta_\square(G_n,\tilde G_n)\to 0$.
{Thus $(\tilde G_n)$ is a sequence of weighted graphs with
nodeweights one which} converges to $U$, implying that it can be
reordered in such a way that $\|W_{\tilde G_n}-U\|_\square\to 0$.
But this means that $G_n/\PP_n$ can be relabeled in such a way that
$\|W_{G_n/\PP_n}-U\|_\square\to 0$, which in turn implies that $G_n$
itself can be {relabeled so} that $\|W_{G_n}-U\|_\square\to 0$, as
desired. }
\end{proof}

The previous lemma suggests that we extend the definition of the distance
$\hat\delta$ to the case when one of the arguments is a graphon $U$:
\begin{equation}
\label{hatdelta-U-G}
\hat\delta_\square(U,G)=\min_{G'}\|U-W_{G'}\|_\square
\end{equation}
where the minimum goes over all relabelings $G'$ of $G$. Then the
lemma asserts that if $G_n$ is a convergent graph sequence with
uniformly bounded edgeweights {and no dominant nodes}, then
$\hat\delta_\square(U,G_n)\to 0$.

{In the special case of nodeweights one}, Theorem \ref{TWIN-SQUARE}
and Lemma \ref{lem:norm-conv-graph} {naively} suggest the stronger
statement that $\hat\delta_\square(U, G)$ can be bounded by a
function $f(\delta_\square(U,W_G),|V(G)|)$ such that $f(x,n)\to 0$ if
$x\to 0$ and $n\to\infty$. However, this is false, as the following
example shows.

\begin{example}\label{ex:norm-metric}
Let $G=K_{n,n}$ be the complete bipartite graph on $2n$ nodes, and
let $H=K_{nm,nm}$ be the complete bipartite graph on $2nm$ nodes,
where $H$ has randomly labeled nodes, and let $U=W_H$. Then
$\delta_\square(U,W_G)=0$. But it is not hard to show that for every
$n$, if $m$ is sufficiently large, then with large probability,
$\|W_H-W_{G'}\|_\square \ge 1/10$ for every relabeling $G'$ of $G$,
implying that $\hat\delta_\square(U,G)\geq 1/10$,
see appendix for details.
\end{example}

\subsection{Convergent Szemer\'edi Partitions}
\label{sec:ConvSzemPart}

Given Lemma~\ref{lem:norm-conv-graph}, we now are ready to
prove Theorem~\ref{thm:Szem-Uniform} about the convergence of
the quotient graphs of suitably chosen Szemer\'edi partitions
for a convergent graph sequence $(G_n)$.

We start by proving the easier direction, namely that
convergent quotients imply convergence of the sequence
$(G_n)$.

Let $\eps>0$, and let $q$ be such that the conditions
(i) and (ii) of Theorem~\ref{thm:Szem-Uniform} hold.  For
a fixed $q$, convergence of $G_n/\PP_n$ is equivalent to convergence
of all edgeweights and nodeweights, which in turn implies convergence in
the $\delta_\square$ distance.  Let $n_0$ be such that
$\delta_\square(G_n/\PP_n,G_m/\PP_m)\leq \eps$ whenever $n,m\geq n_0$,
and $|V(G_n)|\geq q$ whenever $n\geq n_0$.
Then $\delta_\square(G_n,G_m)\leq 3\eps$ for all $n,m\geq n_0$
by the triangle inequality, the property (i) and the fact that
$\delta_\square(G_n,G_n/\PP_n)\leq d_\square(G_n,(G_n)_{\PP_n})$.
This proves that $(G_n)$ is a Cauchy sequence in the metric
$\delta_\square$, and hence left-convergent.

To prove the necessity of the conditions (i) and (ii), consider a
convergent sequence $(G_n)$, a graphon $U'$ such that
$\delta_\square(G_n,U')\to 0$, and a constant $\eps>0$.  With the
help of the Weak Regularity Lemma for graphons,
Lemma~\ref{lem:W-szem}, we can find a partition $\PP'$ of $[0,1]$
into $q_0\leq 2^{10/\eps^2}$ classes such that
$\|U'-U'_{\PP'}\|_\square<  \eps/\sqrt 5$. Applying a
measure-preserving map to both $U'$ and $\PP'$, this allows us to
find a graphon $U$ and a partition $\PP''$ of $[0,1]$ into $q_0$
consecutive intervals such that $\|U-U_{\PP''}\|_\square\leq
\eps/\sqrt 5 $ and $\delta_\square(G_n,U)\to 0$.  Appealing to
Lemma~\ref{lem:norm-conv-graph} we finally relabel the graphs in the
sequence $(G_n)$ in such a way that for the relabeled sequence (which
we again denote by $(G_n)$), we get convergence in norm,
\[\|W_{G_n}-U\|_\square\to 0
\]
as $n\to\infty$.

On the other hand, by {Lemma~\ref{lem:szem-weak}}, we can find a
sequence of partitions $(\PP'_n)$ such that $\PP'_n$ is a weakly
$(\eps/\sqrt 5)$-regular partition of $G_n$ with not more than
$2^{10/\eps^2}$ classes. Let $q_n$ be the number of classes in
$\PP'_n$, and let $q=\max_{n\ge 0} q_n$. By the bound
\eqref{WP-2Opt}, we can refine the partitions $\PP'_n$ to obtain
weakly $(2\eps/\sqrt 5\leq\eps)$-regular partitions $\PP_n''$ with
exactly $q$ classes whenever $|V(G_n)|\geq q$.  In a similar
way, we can refine the partition $\PP''$ to obtain a partition $\PP$
of $[0,1]$ into $q$ consecutive intervals such that
\[
\|U-U_{\PP}\|_\square\leq \frac 2{\sqrt 5}\eps\leq \frac {9}{10}\eps.
\]

Let $n_0$ be such that for $n\ge n_0$
\[
\|W_{G_n}-U\|_\square\leq \frac\eps {30}
\qquad\text{and}\qquad
\frac q{|V(G_n)|}\leq\frac\eps {30},
\]
and let $\alpha_i$ be the Lebesgue measure of the $i^{\text{th}}$
partition class of $\PP$. For $n<n_0$, we then set $\PP_n=\PP''_n$,
and for $n\geq n_0$ we define $\PP_n$ to be the partition into the
classes $V_1^{(n)}=\{1,\dots,k_1^{(n)}\}$,
$V_2^{(n)}=\{k_1^{(n)}+1,\dots,k_1^{(n)}+k_2^{(n)}\}$, etc., where
the integers $k_i^{(n)}$ are chosen in such a way that $\lfloor
\alpha_i |V(G_n)|\rfloor \leq k_i^{(n)}\leq \lceil \alpha_i
|V(G_n)|\rceil$. With this definition, we have that
\[
\|(W_{G_n})_\PP-W_{(G_n)_{\PP_n}}\|_\square
\leq
\|(W_{G_n})_\PP-W_{(G_n)_{\PP_n}}\|_1\leq \frac q{|V(G_n)|}\leq
\frac\eps{30}
\]
for all $n\geq n_0$.  Combined with the triangle inequality and
the bound \eqref{CONTRACT}, this gives
\[
\begin{aligned}
d_\square&(G_n,(G_n)_{\PP_n})
=\|W_{G_n}-W_{(G_n)_{\PP_n}}\|_\square
\\
&\leq
\|W_{G_n}-U\|_\square
+\|U-U_\PP\|_\square
+\|U_\PP-(W_{G_n})_{\PP}\|_\square
+\|(W_{G_n})_{\PP}-W_{(G_n)_{\PP_n}}\|_\square
\\
&\leq
2\|W_{G_n}-U\|_\square
+\|U-U_\PP\|_\square
+\|(W_{G_n})_\PP-W_{(G_n)_{\PP_n}}\|_\square
\leq 2\frac\eps{30} +\frac{9}{10}\eps+\frac\eps{30}=\eps
\end{aligned}
\]
whenever $n\geq n_0$.  Thus $\PP_n$ is a weakly $\eps$-regular
partition of $G_n$, whether $n<n_0$, or $n\geq n_0$. By definition,
it also is a partition into exactly $q$ classes whenever
$|V(G_n)|\geq q$.  This proves (i).

To prove (ii), we note that for $n\geq n_0$, we have
\[
d_\square(G_n/\PP_n,U/\PP)
\leq
\|W_{(G_n)_{\PP_n}}-(W_{G_n})_\PP\|_\square
+\|(W_{G_n})_\PP-U_\PP\|_\square
\leq
\frac q{|V(G_n)|}
+\|W_{G_n}-U\|_\square.
\]
As $n\to 0$, the right hand side goes to $0$, as required.

\section{Parameter Testing}
\label{sec:testing}

\subsection{Definitions and Statements of Results}

In this section we discuss the notion of continuous graph parameters,
and the closely related notion of parameter testing. In parameter
testing, one wants to determine some parameter of a large, {simple}
graph $G$, e.g., the edge density or the density of the maximum cut.
It is usually difficult to determine the exact value of such a
parameter, but using sufficiently large samples, one might hope to
approximate the value of parameter with large probability at a much
lower computation cost; recall Definition (\ref{def:TESTABLE}).

Throughout this section, $G$ will  be a simple graph; this will
be made explicit only in the statements of theorems and supplements.
Our principal theorem gives several useful characterizations of
testable graph parameters.

\begin{theorem}\label{TESTABLE-PAR}
Let $f$ be a bounded simple graph parameter. Then the following are
equivalent:

\medskip

{\rm  (a)} $f$ is testable.

\medskip

{\rm (b)} For every $\eps>0$, there is an integer $k$ such that for
every {simple} graph $G$ on at least $k$ nodes,
\[
|f(G)-\E(f(\RandG(k,G))|\le\eps.
\]

\medskip

{\rm (c)} For every convergent {sequence $(G_n)$ of simple graphs}
with
$|V(G_n)|\to\infty$, the limit of $f(G_n)$ exists as $n\to\infty$.

\medskip

{\rm (d)} There exists a functional $\check{f}(W)$ on $\WW$ that is
continuous in the rectangle norm, and $\check{f}(W_G)-f(G)\to 0$ if
$|V(G)|\to\infty$.

\medskip

{\rm (e)} For every $\eps>0$ there is an $\eps_0>0$ and an
$n_0\in\Z_+$ such that if $G_1,G_2$ are two {simple}
graphs with
$|V(G_1)|,|V(G_2)|\ge n_0$ and $\delta_\square(G_1,G_2)<\eps_0$, then
$|f(G_1)-f(G_2)|<\eps$.
\end{theorem}

If we want to use (e) to prove that a certain invariant is testable,
then the complicated definition of the $\delta_\square$ distance may
cause a difficulty. So it is useful to show that (e) can be replaced
by a weaker condition, which consists of three special cases of (e).
{For this purpose, we define the disjoint union $G\cup G'$ of two
graphs $G$ and $G'$ as the graph whose node set is the disjoint union
of $V(G)$ and $V(G')$, and whose edge set is $E(G)\cup E(G')$.}

\begin{suppl}\label{CONT-INF-SUPP1}
The following three conditions together are also equivalent to
testability:

\medskip

{\rm (e.1)} For every $\eps>0$ there is an $\eps'>0$ such that if $G$
and $G'$ are two {simple} graphs on the same node set and
$d_\square(G,G')\le\eps'$, then $|f(G)-f(G')|<\eps$.

\smallskip

{\rm (e.2)} For every simple graph $G$, $f(G[m])$ has a limit as
$m\to\infty$.

\smallskip

{\rm (e.3)} $f(G{\cup} K_1)-f(G)\to 0$ if $|V(G)|\to\infty$.
\end{suppl}

We formulate two further conditions for testability, in terms of
Szemer\'edi partitions.
Recall that a partition $\{V_1,\dots,V_k\}$
of a finite set $V$ is an {\it equitable partition} if $\lfloor
|V|/k\rfloor \le |V_i|\le \lceil |V|/k\rceil$ for every $1\le i\le
k$. Let $\PP=\{V_1,\dots,V_k\}$ be an equitable partition of
the node set of a simple graph $G$. A pair $(V_i,V_j)$ if partition
classes is called an {\it $\eps$-regular pair}, if for all
$X\subseteq V_i$ and $Y\subseteq V_j$ with $|X|,|Y|\ge\eps
|V(G)|/k$, we have
\[
\Biggl|
\frac{e_G(X,Y)}{|X|\cdot|Y|}-\frac{e_G(V_i,V_j)}{|V_i|\cdot|V_j|}\Biggr|
\le\eps.
\]
The partition $\PP$ is {\it $\eps$-regular} if all but at most
$\eps k^2$ pairs $(V_i,V_j)$ are $\eps$-regular.

Every $\eps$-regular partition is weakly $(7\eps)$-regular,
but in the reverse direction only a much weaker implication
holds: a weakly $\eps$-regular partition with $k$ classes is
$\tilde\eps$-regular with $\tilde\eps =\sqrt[3]{k^2\eps}$.

The ``original'' Szemer\'edi Lemma can be stated as follows.

\begin{lemma}[Szemer\'edi Regularity Lemma
\cite{Szem}]\label{lem:szem-strong}
For every $\eps>0$ and $l>0$ there is a $k(\eps,l)>0$ such that
every simple graph $G=(V,E)$ with at least $l$ nodes has an
$\eps$-regular partition into at least $l$ and at most $k(\eps,l)$
classes.
\end{lemma}

Let $G$ be a graph and let $\PP$ be an equitable partition of
$V(G)$. Then $G/\PP$ is a weighted graph with almost equal
nodeweights. We modify this graph by making all nodeweights equal to
1. This way we get a weighted graph $G\div\PP$ with nodeweights $1$
and edgeweights in $[0,1]$.

For every bounded, simple graph parameter $f$ and every weighted
graph $H$ with nodeweights $1$ and edgeweights in $[0,1]$, define
\[
\iv{f}(H)= \E(f(\bG(H)),
\]
where $\bG(H)$ is the graph obtained by the randomizing procedure
described in Section~\ref{sec:RandWeightG}. Clearly $\iv{f}$ is an
extension of $f$.

\begin{suppl}\label{CONT-INF-SUPP2}
Either one of the following conditions is also equivalent to
testability:

\medskip

{\rm (f)} For every $\eps>0$ there is an $\eps'>0$  and an
$n_0\in\Z_+$ such that if $G$ is any graph with $|V(G)|\ge n_0$ and
$\PP$ is an equitable weakly $\eps'$-regular partition of $G$ with
$n_0\leq |\PP|\leq \eps'|V(G)|$, then
$|f(G)-\iv{f}(G\div\PP)|\le\eps$.

\medskip

{\rm (g)} The parameter $f$ has an extension $\check{f}$ to (finite)
weighted graphs with nodeweights $1$ and edgeweights in $[0,1]$ such
that

\smallskip

{\rm (g.1)} for every fixed $n$, $\check{f}$ is a continuous function
of the edgeweights on $n$-node graphs, and

\smallskip

{\rm (g.2)} for every $\eps>0$ there is an $\eps'>0$  and an
$n_0\in\Z_+$ such that if $G$ is any graph with $|V(G)|\ge n_0$ and
$\PP$ is an $\eps'$-regular partition of $G$ with $n_0\leq |\PP|\leq
\eps'|V(G)|$, then $|f(G)-\check{f}(G\div\PP)|\le\eps$.
\end{suppl}

Condition (g) is {\it a priori} weaker than (f) on two counts:
First, (g) allows an arbitrary extension of $f$ to weighted graphs,
while (f) makes assumptions about a specific extension
$\widehat{f}$. Second, (g) states a condition about $\eps$-regular
partitions, for which the condition may be easier to verify than for
weakly $\eps$-regular partitions. Of course, we could formulate two
``intermediate'' conditions, in which only one of these relaxations
is used.

Condition (g.1) concerns graphs on a fixed finite set, so when we say
that $\check{f}$ should be ``continuous'', we do not have to specify
in which metric we mean this. But to be concrete, we will use the
$d_\square$ distance as a metric on these graphs.

\subsection{Proofs}

Before proving the equivalence of the above conditions, we state and
prove a simple lemma.

\begin{lemma}\label{lem:diff-vert-weights}
Let $G,G'$ be weighted graphs with edgeweights in some interval
$I$ of length $|I|$
and total nodeweight one.
If $G$ and $G'$ have equal
edgeweights but different nodeweights,
then
\[
\delta_\square(G,G')\leq
|I|\sum_i|\alpha_i(G)-\alpha_i(G')|.
\]
\end{lemma}

\begin{proof}
Let $\alpha_i$ be the nodeweights of $G$, and
$\alpha_i'$ be those of $G'$.  Consider a coupling
$X$ such that $X_{ii}=\min\{\alpha_i,\alpha_i'\}$.
Then
\[
\begin{aligned}
\delta_\square(&(G,G')\leq
d_\square(G[X],G'[X\T])
\leq
\sum_{i,j,k,l}
X_{ik}X_{jl}|\beta_{ij}(G)-\beta_{kl}(G')|
\leq |I|
\sum_{\bsa i\neq j\text{ or}\\ k\neq l\esa}X_{ik}X_{jl}
\\
&=|I|\Bigl(1-\Bigl(\sum_i\min\{\alpha_i,\alpha_i'\}\Bigr)^2\Bigr)
=|I|\Bigl(1-\Bigl(1-\frac 12\sum_i|\alpha_i-\alpha_i'|\Bigr)^2\Bigr)
\leq |I|\sum_i|\alpha_i-\alpha_i'|.
\end{aligned}
\]
\vskip -10pt
\end{proof}

After these preparations, we are ready to prove the results of the last
section.

\bigskip
\noindent{\bf Proof of Theorem~\ref{TESTABLE-PAR}} We first prove
that (a), (b), (c) and (e) are equivalent:

\smallskip

(a)$\Rightarrow$(b): The definition of testability is very similar to
condition (b): it says, in this language, that a random set $S_k$ on
$k$ nodes of $G$ as in (b) satisfies
\[
|f(G)-f(\RandG(k,G))| \leq\eps
\]
with large probability. This clearly implies that this difference is
small on average.

\smallskip

(b)$\Rightarrow$(c): Let $(G_n)$ be a convergent sequence with
$|V(G_n)|\to\infty$. Given $\eps>0$, let $k$ be such that for every
graph $G$ on at least $k$ nodes, $|f(G)-\E(f(\RandG(k,G))|\le\eps$.
Since $G_n$ is convergent, $t(F,G_n)$ tends to a limit for all graphs
$F$ on  $k$ nodes, from which we get that $t_\ind(F,G_n)$ tends to a
limit $t_\ind(F)$ for all graphs on  $k$ nodes. This means that $
\Pr(\RandG(k,G_n)=F)\to t_\ind(F), $ and so
\[
\E(f(\RandG(k,G_n)))\to \sum_F t_\ind(F) f(F)=a_k.
\]
As a consequence, we have that for all sufficiently large $n$,
\[
|f(G_n)-a_k|\le|f(G)-\E(f(\RandG(k,G_n))|+\eps\le2\eps,
\]
so $f(G_n)$ oscillates less than $4\eps$ if $n$ is large enough. This
proves that the sequence $(f(G_n))$ is convergent.

\smallskip

(c)$\Rightarrow$(e): Suppose that (e) does not hold for some
$\eps>0$; then there are two sequences $(G_n)$ and $(G_n')$ of graphs
such that $|V(G_n)|,|V(G_n')|\to\infty$, $\delta_\square(G_n,G_n')\to
0$, but $|f(G_n)-f(G_n')|>\eps$. We may assume that both graph
sequences are convergent; but then ${\delta}_\square(G_n,G_n')\to 0$
implies that the merged sequence $(G_1,G_1',G_2,G_2',\dots)$ is also
convergent, so by (c), the sequence of numbers
$(f(G_1),f(G_1'),f(G_2),f(G_2'),\dots)$ is also convergent, a
contradiction.

(e) $\Rightarrow$(a): Suppose (a) does not hold. Then there is an
$\eps>0$ and a sequence $(G_n)$ of graphs with $|V(G_n)|\geq n$ such
that with probability at least $\eps$ we have that
$|f(G_n)-f(\RandG(n,G_n))|>\eps$ for all $n$. Now choose $\eps_0>0$
and an $n_0\in\Z_+$ in such a way that $|f(G_1)-f(G_2)|<\eps$
whenever $G_1,G_2$ are two graphs with $|V(G_1)|,|V(G_2)|\ge n_0$ and
$\delta_\square(G_1,G_2)<\eps_0$ (this is possible by (e)). But by
Theorem \ref{thm:sample2}~(iii), we have
$\delta_\square(G_n,\RandG(k,G_n))\to 0$ in probability, implying
in particular that for $n$ large enough,
$\delta_\square(G_n,\RandG(k,G_n))<\eps_0$ with probability at
least $1-\eps/2$.  By our choice of $\eps_0$, this implies that for
$n$ large enough, $|f(G_n)-f(\RandG(k,G_n))|<\eps$  with
probability at least $1-\eps/2$, a contradiction.

\smallskip

\noindent We continue with proving that (d) is equivalent to (a),
(b), (c) and (e).

(e)$\Rightarrow$(d): For $W\in\WWz$, define
$\check{f}(W)=\lim_{n\to\infty} f(G_n)$, where $(G_n)$ is any
sequence of graphs such that $G_n\to W$ { and $|V(G_n)|\to\infty$ (by
(e), $\check{f}(W)$ does not depend on the choice of the sequence
$G_n$ as long as $G_n\to W$)}. We prove that this functional has the
desired properties.

Let $\eps>0$, and let $\eps'$ and $n_0$ be as given by (e) with
$\eps$ replaced by $\eps/3$.  To prove continuity, we will prove that
$|\check{f}(W)-\check{f}(W')|\le\eps$ whenever
$\|W-W'\|_\square\le\eps'/3$. Considering a graph sequence tending to
$W$, we can choose a simple graph $G$ such that $|V(G)|\ge n_0$,
$\delta_\square(G,W)<\eps'/3$ and $|f(G)-\check{f}(W)|\le\eps/3$;
similarly, we can choose a simple graph $G'$ such that $|V(G')|\ge
n_0$, $\delta_\square(G',W')<\eps'/3$ and
$|f(G')-\check{f}(W')|<\eps/3$. Then
\[
\delta_\square(G,G')\le\delta_\square(G,W)+\delta_\square(W,W')+
\delta_\square(W',G')\le \eps',
\]
and hence by (e), $|f(G)-f(G')|\le\eps/3$. But then
\[
|\check{f}(W)-\check{f}(W')|\le
|\check{f}(W)-f(G)|+|f(G)-f(G')|+|f(G')-\check{f}(W)|\le \eps
\]
as claimed.  Note that our proof shows that, in fact, $\check{f}$ is
continuous in the $\delta_\square$ metric.

The fact that $\check{f}(W_{G_n})-f(G_n)\to 0$ whenever
$|V(G_n)|\to\infty$ can now be easily verified by contradiction.
Indeed, assume that this is not the case.  Using compactness, we may
choose a subsequence such that $G_n\to W$ for some $W$.  But this
implies that $f(G_n)\to\check{f}(W)$, and by the continuity of
$\check f$, also that $\check{f}(W_{G_n})\to \check{f}(W)$, a
contradiction.

\medskip

(d)$\Rightarrow$(c). Consider a convergent graph sequence $(G_n)$
with $|V(G_n)|\to\infty$, and let $W\in\WW$ be its limit. Then
$\delta_\square(W_{G_n},W)\to 0$, and by
Lemma~\ref{lem:norm-conv-graph}, $\|W_{G_n'}-W\|_\square\to 0$ for a
relabeling of $G_n$. By the continuity of $\check{f}$, we have
$\check{f}(W_{G_n'})-\check{f}(W)\to 0$. By assumption,
$f(G_n)-\check{f}(W_{G_n'})=f(G_n)-\check{f}(W_{G_n}) \to 0$, and so
$f(G_n)-\check{f}(W)\to 0$. This proves that $(f(G_n))$ is
convergent. \hfill$\square$

\bigskip

\noindent{\bf  Proof of
Supplement~\ref{CONT-INF-SUPP1}:}

(e)$\Rightarrow$(e.1), (e.2), (e.3): To see that (e.1) is a special
case of (e), it suffices to note that if $G_1$ and $G_2$ are two
different graphs on the same set of $n$ nodes, and
$d_\square(G,G')\le\eps'$, then $n\ge 1/\sqrt{\eps'}$. For (e.2),
note that $\delta_\square(G[m],G[m'])=0$ for all $m,m'$. For (e.3),
it suffices to verify that $\delta_\square(G,GK_1)\le 1/|V(G)|$.

(e.1), (e.2), (e.3)$\Rightarrow$(c): Suppose that (c) does not hold.
Then there exist two graph sequences $(G_n)$ and $(G_n')$ { with
$|V(G_n)|\to\infty$ and $|V(G_n')|\to\infty$} such that $G_n,G_n'\to
W$, $f(G_n)\to a$ and $f(G_n')\to b$ as $n\to\infty$, but $a\neq b$.

By (e.1), there exists an $\eps>0$ such that if $G$ and $G'$ are two
graphs on the same node set, and $d_\square(G,G')\leq\eps$, then
$|f(G)-f(G')|\leq|a-b|/4$. Let $\eps_1={(\eps/32)^{67}}$. By
Lemma~\ref{lem:short-Szemeredi-simple}, for every $n$ there is a
simple graph $H_n$ whose number of nodes $k$ depends only on $\eps$
such that $\delta_\square(G_n, H_n)\leq \eps_1/2$. By selecting an
appropriate subsequence, we may assume that $H_n=H$ is the same graph
for all $n$. Since $(G_n)$ and $(G_n')$ have the same limit, it
follows that $\delta_\square(G_n',H)\le \eps_1$ for all $n$ that are
large enough.

Let us add to each $G_n$ at most $k-1$ isolated nodes so that the
resulting graph $G^*_n$ has $km_n$ nodes for some integer $m_n$. The
$m_n$-fold blow-up $H[m_n]$ of $H$ then satisfies
$
\delta_\square(G^*_n, H[m_n]) \le \eps_1,
$
and so by Theorem~\ref{TWIN-SQUARE}, for a suitable overlay of $G_n$
and $H[m_n]$, we have
$
d_\square(G^*_n, H[m_n]) \le \eps,
$
and so, by the definition of $\eps$,
\[
|f(G^*_n)-f(H[m_n])| \le \frac{|a-b|}{4}.
\]
Using (e.3), we see that $f(G_n^*)-f(G_n)\to 0$, and hence
\[
|f(G_n)-f(H[m_n])| \le \frac{|a-b|}{3}
\]
if $n$ is large enough. Similarly, $H$ has a $m_n'$-node blow-up
$H[m_n']$ such that
\[
|f(G'_n)-f(H[m_n'])| \le \frac{|a-b|}{3}.
\]
But since $H[m_n]$ and $H[m_n']$ are blow-ups of the same graph $H$,
(e.2) implies that $f(H[m_n])-f(H[m_n'])\to 0$, a contradiction.
\hfill$\square$
\medskip

\noindent{\bf Proof of Supplement \ref{CONT-INF-SUPP2}:}
We will assume (without loss of generality) that $|f|\leq 1$.

(e)$\Rightarrow$(f). Choose $\eps_0$ and $n_0$ so that (i) for
$n\geq n_0$ the bound in Lemma~\ref{lem:GH-CLOSE} is {at most
$\eps_0$ with probability at least $1-\eps/4$} and (ii) for two
graphs $G_1$ and $G_2$ with $|V(G_1)|,|V(G_2)|\ge n_0$ and
$\delta_\square(G_1,G_2)\leq 3\eps_0$, we have
$|f(G_1)-f(G_2)|<\eps/4$. Suppose that $|V(G)|\geq n_0$ and
$n_0\leq|\PP|\leq \eps_0 |V(G)|$. By definition of $\eps_0$-regular
partitions and by Lemma~\ref{lem:diff-vert-weights} we have
$\delta_\square(G,G\div\PP)\le 2\eps_0$. Furthermore, by Lemma
\ref{lem:GH-CLOSE}, we have
$\delta_\square(G\div\PP,\bH(G\div\PP))\le \eps_0$ with probability
at least $1-\eps/4$. If this occurs, then
$\delta_\square(G,\bH(G\div\PP))\le 3\eps_0$, and by the choice of
$\eps_0$ and $n_0$, it follows that $|f(G)-f(\bH(G\div\PP))|\le
\eps/2$. Hence
\[
|f(G)-\iv{f}(G)|=\bigl|\E\bigl(f(G)-f(\bH(G\div\PP))\bigr)\bigr| \le
\frac \eps 2(1-\eps/4) + {2}(\eps/4)<\eps.
\]

\medskip

(f)$\Rightarrow$(g). All we have to verify is that for weighted
graphs on a fixed node set (say $[n]$), $\iv{f}$ is a continuous
function of the edgeweights. Consider two weighted graphs $H_1,H_2$
with $V(H_1)=V(H_2)=[n]$ and $d_1(H_1,H_2)\le \eps/n^2$. By Lemma
\ref{lem:GH-COUPL}, the randomized simple graphs $\bG(H_1)$ and
$\bG(H_2)$ can be coupled so that $\E [d_1(\bG(H_1),\bG(H_2))]\le
\eps/n^2$, which by Markov's inequality implies that
$\bG(H_1)=\bG(H_2)$ with probability at least $1-\eps$. This implies
that $|\iv{f}(\bG(H_1))-\iv{f}(\bG(H_2))|\le \eps$. (Note that this
argument did not require that $f$ is testable.)

\medskip

(g)$\Rightarrow$(e.1),(e.2),(e.3). Let $\eps>0$, and consider the
extension $\check{f}$ postulated in (g). By property (g.2), we can
choose an $\eps_1>0$ and an $n_1\in\Z_+$ such that if $|V(G)|\ge n_1$
and $\PP$ is an $\eps_1$-regular partition of $G$ with $n_1\leq
|\PP|\leq\eps_1 |V(G)|$, then $|f(G)-\check{f}(G\div\PP)|\le\eps/3$.
Using the Regularity Lemma \ref{lem:szem-strong}, fix $k=
k(\eps_1/2,n_1)>n_1$ so that every graph $G$ with at least $n_1$
nodes has an $(\eps_1/2)$-regular partition with at least $n_1$ and
at most $k$ classes. Let $n_0=\max\{n_1,\lceil k/\eps_1\rceil\}$.
Using condition (g.1), choose an $\eps_2>0$ such that for any two
graphs $H_1$ and $H_2$ on the same set of $m\le n_0$ nodes with
$d_\square(H_1,H_2)\le\eps_2$ we have
$|\check{f}(H_1)-\check{f}(H_2)|\le\eps/3$. Finally, choose
$\eps'=\min\{\eps_2/4, \eps_1^3/(8k^2)\}$.

To prove (e.1), let $G$ and $G'$ be two simple graphs on the same
node set with $d_\square(G,G')\le\eps'$. If $n=|V(G)|\leq n_0$, then
$|f(G_1)-f(G_2)|=|\check{f}(G_1)-\check{f}(G_2)|\le\eps/3$, so we can
assume $|V(G)|>n_0$. Let $\PP=(V_1,\dots,V_l)$ be an
$(\eps_1/2)$-regular partition of $G$ and $n_1\leq l\leq k$. From the
assumption that $d_\square(G,G')\le\eps'\le \eps_1^3/(8k^2)$, it
follows that $\PP$ is an $\eps_1$-regular partition of $G'$. By the
choice of $\eps_1$ and $n_1$ it follows that
\[
|f(G)-\check{f}(G\div\PP)|\le\frac{\eps}{3}\qquad\text{and}\qquad
|f(G')-\check{f}(G'\div\PP)|\le\frac{\eps}{3}.
\]
Furthermore, we have
\[
n^2 d_\square(G,G')\ge l^2 d_\square(G\div\PP,G'\div\PP)
\Bigl\lfloor \frac nl \Bigr\rfloor^2,
\]
and so
\[
d_\square(G\div\PP,G'\div\PP) \le \frac{n^2}{l^2}
\Bigl\lfloor \frac nl \Bigr\rfloor^{-2}
d_\square(G,G')\le 4\eps'\leq \eps_2.
\]
Hence by the choice of $\eps_2$, we have
\[
|\check{f}(G\div\PP)-\check{f}(G'\div\PP)|\le\frac{\eps}{3}.
\]
Summing up,
\[
|f(G)-f(G')|\le
|f(G)-\check{f}(G\div\PP)|+|\check{f}(G\div\PP)-\check{f}(G'\div\PP)|
+|\check{f}(G'\div\PP)-f(G')| \le \eps.
\]

The proof of (e.2) is similar, but divisibility concerns cause some
complications. We have to show that $f(G[q])$ is a Cauchy sequence,
i.e., we have to show that given $G$ and $\eps>0$, we can find a
$q_0$ such that $|f(G[q])-f(G[q'])|\le 2\eps$ whenever $q,q'\geq
q_0$. Let $n=|V(G)|$, $p>1$ large enough, and let
$\PP=\{V_1,\dots,V_l\}$ be an $(\eps_1/2)$-regular partition of
$G[p]$ with $n_1\le l\le k$. By our choice of $p$ we have, in
particular, that $l\leq \frac{\eps_1}{2}np$, so we may apply (g.2)
to get
\begin{equation}\label{FGP}
|f(G[p])-\check{f}(G[p]\div\PP)|\le\frac{\eps}{3}.
\end{equation}

Let $r$ be sufficiently large, and let $q=pr+s$ with $0\leq s<p$.
The graph $G[pr]$ arises from $G[p]$ by blowing up each node into
$r$ nodes, and then $G[q]$ arises from $G[pr]$ adding $sn$ further
nodes.

First we consider the graph $G[pr]$. The partition $\PP$ determines
a partition $\QQ=\{U_1,\dots,U_l\}$ of $V(G[pr])$. Let
\[
p_{ij}=\frac{e_{G[p]}(V_i,V_j)}{|V_i|\cdot|V_j|}
=\frac{e_{G[pr]}(U_i,U_j)}{|U_i|\cdot|U_j|}
\]
denote the edge density between $V_i$ and $V_j$ in $G[p]$ (which is
the same as the edge density between $U_i$ and $U_j$ in $G[pr]$).

We claim that if $(V_i,V_j)$ is an $(\eps_1/2)$-regular pair in
$\PP$, then for all $X\subseteq U_i$ and $Y\subseteq U_j$ with
$|X|,|Y|\ge 2\eps_1npr/(3l)$, we have
\begin{equation}\label{PR-REG}
p_{ij}-\frac{\eps_1}2\le \frac{e_{G[pr]}(X,Y)}{|X|\cdot|Y|}\le
p_{ij}+\frac{\eps_1}2.
\end{equation}
For $u\in V_i$, let $x'_u$ denote the number of elements in $X$
among the $r$ copies of $u$ in $G[pr]$,and let $x_u=x'_u/r$. Clearly
$0\le x_u\le 1$ and $\sum_u x_u \ge \lfloor 2\eps_1npr/(3l)\rfloor$.
We define $y_v$ for $v\in V_j$ analogously. Then we have
\begin{equation}\label{SIRREG3}
e_{G[pr]}(X,Y)-|X|\cdot|Y|(p_{ij}-\eps_1) = \sum_{u\in
V_i}\sum_{v\in V_j}r^2 x_uy_v (a_{uv}- p_{ij}-\frac{\eps_1}2).
\end{equation}
Since the right hand side is linear in each $x_u$, it attains its
minimum over $0\le x_u\le 1$, $\sum_u x_u \ge \lfloor
2\eps_1npr/(3l)\rfloor$, at a vertex of this domain, which is a 0-1
vector. Similarly, the minimizing choice of the $y_v$ is a 0-1
vector. Let $S=\{u:~x_u=1\}$ and $T=\{v:~y_v=1\}$, then $|S|,|T|\ge
\lfloor 2\eps_1npr/(3l)\rfloor \ge \frac{\eps_1}2|V(G[p]|/l$. The
right hand side of \eqref{SIRREG3} is equal to
\[
r^2\Bigl(e_{G[p]}(S,T)-|S|\cdot|T|\bigl(p_{ij}
-\frac{\eps_1}2\bigr)\Bigr)\ge0,
\]
since $(V_i,V_j)$ is $(\eps_1/2)$-regular. This proves the first
inequality in \eqref{PR-REG}; the proof of the second is similar.

Now $\QQ$ is not necessarily an equitable partition; the largest and
smallest class sizes may differ by $r$, not by $1$. Let us remove
$r$ nodes from those classes of $\QQ$ that are too big. To get a
partition of $V(G[q])$, we have to add back these nodes ($t\le l$)
and $sn$ further nodes. Let us distribute these nodes as equally as
possible between the classes, to get an equitable partition
$\QQ'=\{U_1',\dots,U_l'\}$ of $G[q]$.

We claim that $\QQ'$ is $\eps_1$-regular. Consider a pair
$(V_i,V_j)$ that is $(\eps_1/2)$-regular in $G[p]$, and subsets
$X'\subseteq U_i'$ and $Y'\subseteq U_j'$ with
$|X|,|Y|\ge\eps_1nq/l$. Remove the new nodes from $X'$ and $Y'$ to
get $X\subseteq U_i$ and $Y\subseteq U_j$. Clearly $|X|,|Y|\ge
\eps_1nq/l - r- \lceil sn/l\rceil\ge 2\eps_1npr/(3l)$, and so
\eqref{PR-REG} is satisfied. Now it is easy to check that
\[
\left|\frac{e_{G[q]}(U'_i,U'_j)}{|U'_i|\cdot|U'_j|}
-\frac{e_{G[pr]}(U_i,U_j)}{|U_i|\cdot|U_j|}\right| \le
\frac{\eps_1}4
\qquad\text{and}\qquad
\left|\frac{e_{G[q]}(X',Y')}{|X'|\cdot|Y'|}
-\frac{e_{G[pr]}(X,Y)}{|X|\cdot|Y|}\right| \le \frac{\eps_1}4
\]
if $p>32 l/\eps_1^2$ and $q>32 p/\eps_1^2$, which proves that
\[
\left|\frac{e_{G[q]}(U'_i,U'_j)}{|U'_i|\cdot|U'_j|}
-e_{G[q]}(X',Y'){|X'|\cdot|Y'|}\right| \le \eps_1.
\]

Thus $\QQ'$ is an $\eps_1$-regular  partition of $G[q]$. It follows
by (g.2) and the choice of $\eps_1$ that
\begin{equation}\label{FGP-2}
|f(G[q])-\check{f}(G[q]\div\QQ')|\le \frac{\eps}{3}.
\end{equation}
The weighted graphs $G[p]\div\PP$ and $G[q]\div\QQ'$ are very close
to each other. In fact, $G[p]\div\PP\cong G[pr]\div\QQ$, while it is
easy to check that $d_\square(G[pr]\div\PP, G[q]\div\QQ)<\eps_2$ if
$q>4p/\eps_2$, and so (g.1) implies that
\begin{equation}\label{FGP-3}
|\check{f}(G[p]\div\PP)-\check{f}(G[q]\div\QQ')|\le\frac{\eps}{3}.
\end{equation}
Now \eqref{FGP}, \eqref{FGP-2} and \eqref{FGP-3} imply that
$|f(G[p])-f(G[q])|\le\eps$. So if
$q,q'>\max(32p/\eps_1^2,4p/\eps_2\}$, then
\[
|f(G[q])-f(G[q'])|\le |f(G[q])-f(G[p])|+|f(G[p])-f(G[q'])|\le 2\eps.
\]
Thus $f(G[q])$ is a Cauchy sequence, which proves (e.2).

The proof of (e.3) is similar but easier, and is left to the
reader.  \hfill$\square$

\section{Concluding Remarks}
\label{sec:remarks}

\subsection{Norms Related to the Cut-Norm}\label{MORE-NORMS}

The cut-norm (and the cut-distance of graphs) is closely related to
several other norms that are often used. We formulate these
connections for the case of graphons, but similar remarks would
apply to the $d_\square$ distance of graphs.

We start with the remark that we could restrict the sets $S,T$ in
Definition \ref{CUTNORM-W}; it turns out that ``reasonable''
restrictions only change the supremum value by a constant factor. In
particular, it is not hard to see that
\begin{equation}\label{SQUARE-SS}
\frac{1}{2}\|W\|_\square
\le \sup_{S=T}\biggl|\int_{S\times T}
W\biggr| \le \|W\|_\square,
\end{equation}
\begin{equation}\label{SQUARE-ST-DISJ}
\frac{1}{4} \|W\|_\square
\le \sup_{S\cap T=\emptyset}\biggl|\int_{S\times T}
W\biggr|
\le \|W\|_\square,
\end{equation}
and
\begin{equation}\label{SQUARE-ST-COMPL}
{\frac{2}{3} \sup_{S\cap T=\emptyset}
\biggl|\int_{S\times T} W\biggr|
\le\sup_{S=[0,1]\setminus T}\biggl|\int_{S\times T} W\biggr|
\le \sup_{S\cap T=\emptyset}\biggl|\int_{S\times T} W\biggr|,}
\end{equation}
see appendix for  details.

To relate the cut norm to homomorphisms, we start with noticing that
\begin{equation}\label{C4-L2}
t(C_4,U)^{1/4}=(\text{Tr}\ T_U^4)^{1/4},
\end{equation}
and hence the functional $t(C_4,.)^{1/4}$ defines a norm, called the
{\it trace norm} or {\it Schatten norm}, on $\WW$. (More generally,
any even cycle leads to a norm in this way.) The following lemma
shows that the norm in \eqref{C4-L2} is intimately related to the
cut-norm.

\begin{lemma}\label{lem:C4-SQUARE}
For $U\in\WW$ with $\|U\|_\infty\le 1$, we have
\[
\frac14 t(C_4,U)\le \|U\|_\square\le t(C_4,U)^{1/4}.
\]
\end{lemma}

\begin{proof}
The first inequality is a special case of \eqref{LEFTNORM}. To prove
the second inequality, we use (\ref{CUT-IO-W}) and (\ref{C4-L2}):
\[
\|U\|_\square \le \|U\|_\IO = \sup_{|f|,|g|\le 1} \langle f,
T_Ug\rangle,
\]
and here
\begin{align*}
\langle f, T_Ug\rangle &\le \|f\|_2\cdot\|T_Ug\|_2 =
\|f\|_2\cdot\langle T_Ug,T_Ug\rangle^{1/2} = \|f\|_2\cdot\langle g,
T_U^2g\rangle^{1/2}\\
&\le \|f\|_2\cdot\|g\|_2\cdot \|T_U^2\|_{2\to 2}^{1/2}\le
\|T_U^2\|_2^{1/2}=(\text{Tr}\ T_U^4)^{1/4}= t(C_4,U)^{1/4}.
\end{align*}
\end{proof}

Lemma~\ref{lem:C4-SQUARE} allows for a significant simplification of
the proof that the sample from a weighted graphs is close to the
original graph. Indeed, it is possible to use the following easy
lemma instead of the quite difficult Theorem~\ref{thm:NormSample}
(or the equally difficult results of \cite{AFKK}) to establish a
slightly weaker version of Theorem~\ref{thm:sample2}, which is still
strong enough to prove the equivalence of left-convergence and
convergence in metric (see \cite{BCLSV-1} for details).

\begin{lemma}\label{lem:Norm-Sampl}
Let $0<\eps<1$, $0<\delta<1$, and let
$U\in\WW$, $\|U\|_\infty\le 1$.
If
\smallskip
\[
\|U\|_\square\leq \frac 18\eps^{1/4}
\quad\text{and}\quad k\ge \frac{352}{\eps^8}\ln(\frac{2}{\delta}),
\]
then
\[
\Prob\Bigl(\|\bH(k,U)\|_\square \le \eps\Bigr)\geq 1-\delta.
\]
\end{lemma}

\begin{proof}
Applying first Lemma \ref{lem:C4-SQUARE} and then Lemma \ref{lem:t-conc}, we
have
\[
\|\bH(k,U)\|_\square \le t(C_4, \bH(k,U))^{1/4}
\leq \Bigl(t(C_4,U)+\frac{\eps^4}2\Bigr)^{1/4}
\]
with probability at least $1-\delta$.
Thus, using Lemma
\ref{lem:C4-SQUARE} again, we get
\[
\|\bH(k,U)\|_\square
\le\Bigl(t(C_4,U)+\frac{\eps^4}2\Bigr)^{1/4}\le
\Bigl(4\|U\|_\square+\frac{\eps^4}2\Bigr)^{1/4}\leq \eps,
\]
as claimed.
\end{proof}

A substantial advantage of the norm $t(C_4,W)$ over the cut-norm is
that for weighted graphs (and for many other types of graphons $W$,
for example, for polynomials), it is polynomial-time computable. A
better polynomial-time computable approximation is the {\it
Grothendieck norm}, which approximates the cut-norm within a
constant factor (see \cite{AN}).

\subsection{A Common Generalization of Lemma~\ref{lem:szem-weak}
and Lemma~\ref{lem:W-szem}} \label{sec:WRL-common}

As stated, Lemma~\ref{lem:W-szem} does not generalize
Lemma~\ref{lem:szem-weak} since the partition in
Lemma~\ref{lem:W-szem} is not necessarily aligned with the steps of
the stepfunction $W_G$.  But the following lemma implies both
Lemma~\ref{lem:szem-weak} and Lemma~\ref{lem:W-szem}.

\begin{lemma}\label{lem:SZEM-ALG}
Let $\AA$ be an algebra of measurable subsets of $[0,1]$, and let
\[
\|W\|_{\AA,\square}= {\sup_{S,T\in\AA}}\Bigl|\int_{S\times T}
W\Bigr|.
\]
Then for every graphon $W$ and every $\eps>0$, there exists a
partition $\PP$ of $[0,1]$ into sets in $\AA$ with at most $4^{\lceil
1/\eps^{2}\rceil -1}$ classes such that
\[
\|W-W_\PP\|_{\AA,\square}\le\eps\|W\|_2.
\]
\end{lemma}

\begin{proof}
{Let $\PP$ of be a partition of $[0,1]$ into $q$ classes in $\AA$,
let $S,T\in\AA$, and let $\PP'$ be the partition generated by $S,T$
and $\PP$. Clearly $\PP'$ has at most $4q$ classes, all of which lie
in $\AA$. Since $W_{\PP'}$ gives the best $L_2$-approximation of $W$
among all step functions with steps $\PP'$, we conclude that for
every real number $t$, we have
\[
\|W-W_{\PP'}\|_2^2 \le \|W-W_\PP-t\one_{S\times T}\|^2_2.
\]
Bounding the right hand side by $\|W-W_\PP\|_2^2 - 2t
\langle\one_{S\times T}, W-W_\PP\rangle +t^2$ and choosing
$t=\langle\one_{S\times T}, W-W_\PP\rangle$, this gives
\[
\langle\one_{S\times T}, W-W_\PP\rangle^2 \le
\|W-W_\PP\|_2^2-\|W-W_{\PP'}\|_2^2 =\|W_{\PP'}\|_2^2-\|W_{\PP}\|_2^2
\]
Taking the supremum over all sets $S,T\in\AA$, this gives
\begin{equation}
\label{S-disguise} \|W-W_\PP\|_{\AA,\square}^2 \le
\sup_{\PP'}\|W_{\PP'}\|_2^2- \|W_\PP\|_2^2,
\end{equation}
where the supremum goes over all partitions of $[0,1]$ into $4q$
classes in $\AA$. From this bound, the lemma then follows by standard
arguments.}
\end{proof}

\subsection{Right Convergence}

When studying homomorphisms from $G$ into a small graph $H$ it will
be convenient to consider graphs $H$ with nodeweights,
$\alpha_i(H)>0$, and edgeweights $\beta_{ij}(H)\in\R$ (with $i$
running over all nodes in $H$, and $ij$ running over all edges of
$F$). For such a graph, we define
\[
\hom(G,H)=\sum_{\phi} \prod_{i\in V(G)}\alpha_{\phi(i)}(H)
\prod_{ij\in E(G)}\beta_{\phi(i),\phi(j)}(H)
\]
where the sum runs over all maps $\phi$ from $V(G)$ to $V(H)$ and
$\beta_{\phi\phi'}(H)$ is set to zero if $\phi\phi'$ is not an edge
in $H$. We call $\hom(G,H)$ the weighted number of $H$-colorings of
$G$.

For dense graphs $G$, the weighted $H$-coloring numbers $\hom(G,H)$
turn out to be most interesting when all edgeweights of $H$ are
strictly positive (we say that $H$ is a soft-core graph if this is
the case).  Under this assumption, the homomorphism numbers
$\hom(G,H)$ typically grow exponentially in the number of edges in
$G$. For homomorphism into small graphs $H$, it is therefore natural
to consider the logarithm of $\hom(G,H)$ divided by the number of
nodes in $G$ to the power two. We will also consider the
``microcanonical ensemble'', where the number of nodes in $V(G)$
colored by a given color $c\in V(H)$ is fixed to be some constant
proportion $a_c$ of all nodes. We denote these homomorphism densities
by $\hom_{\ba}(G,H)$, where $\ba$ is the  vector with components
$a_c$, $c\in V(H)$.

The above perspective leads to two {\it a priori} different notions
of convergence: A sequence of weighted graphs $(G_n)$ will be called
{\it convergent from the left} if the homomorphism densities
$t(F,G_n)$ converge for all finite graphs $F$, and it will be called
{\it convergent from the right} if the quantity
$|V(G_n)|^{-2}\log\hom_{\ba}(G_n,H)$ converges {for all $\ba$ and}
all soft-core graphs $H$. It will turn out that these two notions
are closely related. Convergence from the left implies convergence
from the right (both for the standard homomorphism numbers and the
microcanonical ones), and convergence from the right for the
microcanonical homomorphism numbers implies convergence from the
left.  This will be discussed in more detail in the continuation of
this paper \cite{dense2}.

\newpage
\section{Appendix}
\label{sec:appendix}

In this appendix, we collect various proofs which
were omitted in the body of the paper.

\subsection{Proof of Corollary~\ref{cor:unif-Sem}}

{\bf Proof of (i)}:
Let $\PP'$ be a partition of $[0,1]$ into $q'\leq 2^{81/(8\eps^{2})}$
classes such that $\|W-W_{\PP'}\|_\square\leq \frac {4\eps}9
\|W\|_2$. Then there exists an equipartition $\PP$ of $[0,1]$ into
$q$ classes such that at most $q'$ of its classes intersect more than
one class of $\PP'$. Let $R$ be the union of these exceptional
classes, and let $U$ be the step function which is equal to
$W_{\PP'}$ on $([0,1]\setminus R)^2$, and $0$ on the complement. Then
$\|W-U\|_\square$ can easily be bounded by decomposing the sets $S,T$
in the definition of the cut-norm into the part contained in
$[0,1]\setminus R$ and its complement.  Using the fact that
$\lambda(R)\leq \frac {q'}q\leq 2^{-79/(8\eps^2)}\leq \eps^2
2^{-79/8}$, this leads to the estimate
\[
\|W-U\|_\square\leq \Bigl(\frac {4\eps}9 +
\sqrt{2\lambda(R)}\Bigr) \|W\|_2 \leq
\frac\eps 2 \Bigl(\frac 89 + \sqrt{8\cdot 2^{-79/8}}\Bigr)\|W\|_2
\leq \frac\eps 2 \|W\|_2.
\]
By construction, $U$ is a step function with steps in $\PP$. Using
the bound \eqref{WP-2Opt}, we conclude that $\|W-W_\PP\|_\square\leq
2\|W-U\|_\square\leq\eps\|W\|_2$, which gives the first statement of
the corollary.
The second assertion in statement (i) is proved in a
similar way, starting from a common refinement of $\PP'$ and
$\tilde\PP$.

{\bf Proof of  (ii)}: Let $V=V(G)$, let $n=|V|$, and assume
without loss of generality that $\alpha_G=1$. Choosing a partition
$\PP'=(V_1',\dots,V_{q'}')$ of $V$ with $q'\leq 2^{81/(8\eps^2)}$
classes such that $d_\square(G,G_{\PP'})\leq \frac {4\eps}9\|G\|_{2}
$, we would like to divide each class in $\PP'$ into subclasses
$V_i$ such that all of them obey the condition \eqref{weight-cond}.
To this end, we proceed as follows:  Starting with $V_1'$, we
successively remove sets $V_1,V_2,\dots$ from first $V_1'$, then
$V_2'$, etc. so that
\begin{equation}
\label{size-const}
\Bigl|\alpha_{G[V_i]}-\frac {1}q\Bigr|
<\alpha_{\max}(G)
\qquad\text{for all $i=1,\dots,q$}
\end{equation}
and after each step
\[
\Bigl|\sum_{i=1}^t \alpha_{G[V_i]}
-\frac{t}q\Bigr|\leq \frac{\alpha_{\max}(G)}2.
\]
When it is not possible to further remove a set $V_i$ from $V_1'$
while maintaining these constraints, we are left with a (possibly
empty) remainder $R_1$ that has weight $\alpha_{G[R_1]}<1/q$
(otherwise, we could have continued for at least one more step).
Continuing with $V_2',\dots,V_{q'}'$, we will eventually end up with
$q-q'\leq t\leq q$ disjoint sets $V_1,\dots, V_t$ obeying the
condition \eqref{size-const}, and $r\leq q'$ non-empty remainders
$R_i$ such that the total weight of their union, $R=\bigcup_iR_i$,
obeys the condition $|\alpha_{G[R]}-(q-t)/q|\leq \alpha_{\max}/2$.
Using this condition, it is not hard to see that $R$ can be split
into $q-t$ final sets $V_{t+1},\dots, V_q$ obeying the condition
\eqref{size-const}. Since each of the remainders  had weight at most
$1/q$, the total weight of $R$ is at most $q'/q\leq
2^{-79/(8\eps^2)}\leq \eps^2 2^{-79/8}$.

From here on the proof is completely analogous to the  proof of
(i).  \hfill$\square$

\subsection{Proof of Lemma~\ref{lem:short-Szemeredi-simple}}

We start with the observation that the proof of the last
section actually gives the stronger bound
\[
\|W-W_\PP\|_\square\leq
\eps\Bigl(\frac 89 + \sqrt{8\cdot 2^{-79/8}}\Bigr)\|W\|_2
\leq 0.982\eps\|W\|.
\]
Applying this bound to $W_G$, this gives a weighted
graph $H=W_G/\PP$ on $q\geq 2^{20/\eps^2}\geq 2^{20}/\eps^2$ nodes
such that $H$ has nodeweights one and
$\delta_\square(G,H)\leq 0.982\eps$.
Combined with Lemma~\ref{lem:GH-CLOSE}, this gives the existence of a weighted
graph $\widetilde H$ on $[q]$ such that

\[
\delta_\square(G,\widetilde H)\leq 0.982\eps+
\frac{4}{\sqrt q}\leq \eps,
\]
as required.  \hfill$\square$

\subsection{Proof of Lemma~\ref{lem:delta=delta}}

We notice that $\delta_\square(U,W)$ as well as the infima and limits
in (\ref{eq:DD1}), (\ref{eq:DD2}) and (\ref{eq:DD3}) are continuous
in both $U$ and $W$, with respect to the $\|.\|_\square$ norm. This
fact and Lemma~\ref{lem:as-approx} imply that it is enough to prove
the lemma for graphons $U$ and $W$ that are interval step functions
with equal steps, corresponding to some finite graphs $G$ and $G'$.

Furthermore, the inequalities
\[
\inf_{\phi,\psi}\|U^\phi-W^\psi\|_\square \le
\inf_{\psi}\|U-W^\psi\|_\square\le
\liminf_{n\to\infty}\min_{\pi}\|U-W^{\tilde{\pi}}\|_\square
\]
are trivial, so it suffices to prove that
\begin{equation}\label{eq:DDpf1}
\delta_\square(U,W)\le\|U^\phi-W^\psi\|_\square
\end{equation}
for all measure-preserving maps $\phi,\psi:~[0,1]\to[0,1]$, and
\begin{equation}\label{eq:DDpf2}
\limsup_{n\to\infty}\min_{\pi}\|U-W^{\tilde{\pi}}\|_\square
\le\delta_\square(U,W),
\end{equation}
where the infimum is over all permutations of $[n]$.

To prove \eqref{eq:DDpf1}, we consider two weighted graphs $G$ and
$G'$ with $\alpha_G=\alpha_G'=1$, together with the step functions
$U=W_G$ and $W=W_{G'}$. Let $I_1,\dots, I_n$ be the intervals
$[0,\alpha_1(G)]$, $(\alpha_1(G),\alpha_1(G)+\alpha_2(G)]$, $\dots$,
$(\alpha_1(G)+\dots+\alpha_{n-1}(G),1]$, and similarly for
$I_1',\dots,I_{n'}'$. For two measure-preserving maps
$\phi,\psi:~[0,1]\to[0,1]$, we then rewrite the norm on the right
hand side of \eqref{eq:DDpf1} in the explicit form
\[
\|U^\phi-W^\psi\|_\square = \sup_{S,T\subset
[0,1]}\Bigl|\int_{S\times T}
\bigl(W_G(\phi(x),\phi(y))-W_{G'}(\psi(x),\psi(y))\bigr)\,dx\,dy\Bigr|.
\]
The supremum on the right hand side is attained when $S$ and $T$ are
unions of sets of the form $V_{iu}=
\phi^{-1}(I_i)\cap\psi^{-1}(I'_u)$, $i\in V(G)$, $u\in V(G')$. But
for these sets, the integral on the right is a sum of terms of the
form
\[
\int_{V_{iu}\times V_{jv}}
\bigl(W_G(\phi(x),\phi(y))-W_{G'}(\psi(x),\psi(y))\bigr)\,dx\,dy =
\beta_{ij}(G)\beta_{uv}(G')X_{iu}X_{jv}
\]
where $X_{iu}$ is the Lebesgue measure of the set $V_{iu}$.  As a
consequence, we have that
\[
\|U^\phi-W^\psi\|_\square = d_\square(G[X],G'[X\T]).
\]
Using the fact that $\phi$ and $\psi$ are measure-preserving, it is
not hard to check that $X$ is a coupling of the distributions
$(\alpha_i(G))_{i\in V(G)}$ and $(\alpha_u(G'))_{u\in V(G')}$,
implying that
\begin{align}
\|U^\phi-W^\psi\|_\square = d_\square(G[X],G'[X\T]) \ge
\delta_\square(U,W).
\end{align}

To show \eqref{eq:DDpf2}, let us consider a coupling $X\in
\XX(G,G')$. First construct a special measure-preserving bijection
$\psi:~[0,1]\to[0,1]$ such that
$\|W_G-W_{G'}^\psi\|_\square=d_\square(G[X],G'[X\T])$. As shown
above, this is equivalent to finding a measure-preserving map such
that $X_{iu}$ is the Lebesgue measure of $I_i\cap\psi^{-1}(I'_u)$.
But the construction of such a map is straightforward. Indeed, for
$iu\in V(G)\times V(G')$, let
$b_i=\alpha_1(G)+\dots+\alpha_{i-1}(G)$,
$b'_u=\alpha_1(G')+\dots,\alpha_{u-1}(G')$,
$c_{iu}=b_i+X_{i1}+X_{i2}+\dots+X_{i(u-1)}$ and
$c_{iu}'=b'_i+X_{1u}+X_{2u}+\dots+X_{nu}$. Let $I_{iu}$ and $I'_{iu}$
be the intervals $I_{iu}=(c_{iu},c_{iu}+X_{iu}]$ and
$I'_{iu}=(c'_{iu},c'_{iu}+X_{iu}]$. We then choose $\psi$ to be the
translation that maps $I_{iu}$ into $I'_{iu}$. Then
$\psi^{-1}(I'_u)=\bigcup_i\psi^{-1}(I'_{iu})=\bigcup_iI_{iu}$ and
$I_i\cap \psi^{-1}(I'_u)=I_{iu}$, implying in particular that this
set has measure $X_{iu}$, as required.

So we have two partitions $\{I_1,\dots,I_m\}$ and $\{J_1,\dots,J_m\}$
of $[0,1]$ into intervals, and $\psi$ maps each $I_k$ onto a
$J_{f(k)}$ by translation. Furthermore, we also know that both $U$
and $W$ are constant on each rectangle $I_k\times I_l$ as well as on
$J_k\times J_l$.

Let $N$ be a large integer, and consider the partition
$\{L_1,\dots,L_N\}$ of $[0,1]$ into intervals of size $1/N$. We
define a permutation $\pi$ of $[N]$. For every $k\le m$, the
intervals $I_k$ and $J_{f(k)}$ have the same length, and so the
numbers of intervals $L_i$ contained in them can differ by at most
one. Let $\pi$ match the indices $i$ such that $L_i\subseteq I_k$
with the indices $j$ such that $L_j\subseteq J_{f(k)}$, with at most
one exception. This way $\pi(i)$ is defined for at least $N-3m$
integers $i\in[N]$. Call the corresponding intervals $K_i$ {\it
well-matched}. We extend $\pi$ to a permutation of $[N]$ arbitrarily.

We see that $W^{\psi}(x,y)=W^{\tilde{\pi}}(x,y)$ whenever both $x$
and $y$ belong to well-matched intervals. Hence
\[
\|W^{\psi}-W^{\tilde{\pi}}\|_\square\le
\|W^{\psi}-W^{\tilde{\pi}}\|_1 \le \frac{6m}{N}\|W\|_\infty,
\]
and so
\[
\|U-W^{\tilde{\pi}}\|_\square \le \|U-W^\psi\|_\square +
\|W^{\psi}-W^{\tilde{\pi}}\|_\square \le
\delta_\square(U,W)+\frac{6m}{N}\|W\|_\infty.
\]
This implies \eqref{eq:DDpf2}.  \hfill$\square$

\subsection{Proof of Lemma~\ref{lem:t-conc}}

Starting with the proof of \eqref{WH-Concentr}, let us assume
that $k^2/n\leq \eps/(11\log 2)$ (otherwise the bound
\eqref{WH-Concentr} is trivial).  Using the bound
\eqref{random-W-Exp-t}, we then estimate the probability on the left
hand side of \eqref{WH-Concentr} by
\[
\begin{aligned}
&\Prob\Bigl(\Bigl|t(F,\bH(n,W))-\E[t(F,\bH(n,W))]\Bigr|
>\Bigl(\eps-\frac{2k^2}n\Bigr)
\Bigr)
\\
&\qquad\leq \Prob\Bigl(\Bigl|t(F,\bH(n,W))-\E[t(F,\bH(n,W))]\Bigr|
>\eps\Bigl(1-\frac
1{11\log 2}\Bigr)\Bigr).
\end{aligned}
\]
Applying Lemma~\ref{lem:Azuma} (i) and the observation that
$t(F,W[\{z_1,\dots,z_n\}])
$ changes by at most $2k/n$ if we change one of the variables
$z_i$
we immediately obtain the bound \eqref{WH-Concentr}.

Expressing $\bG(n,W)$ as a function of the random variables
$Z_1,\dots,Z_n$ introduced above, and observing that $t(F,\bG(n,W))$
changes by at most $k/n$ if we change one of the variables $Z_i$, the
proof of \eqref{WG-Concentr} is virtually identical to that of
\eqref{WH-Concentr}. We leave the details to the reader.
\hfill$\square$

\subsection{Detail Concerning Example~\ref{ex:norm-metric}}

Let $H=K_{nm,nm}$ and $G=K_{n,n}$ be the graphs considered in
Example~\ref{ex:norm-metric}.  Here we show that for  every $n$,
there exists an $m$ such that with large probability,
$\|W_H-W_{G'}\|_\square \ge 1/10$ for every relabeling $G'$ of $G$.

Indeed, let $S'$ and $T'$ be
the color classes of $G'$, and let $I$ and $J$ be the subsets of
$[0,1]$ corresponding to $S'$ and $T'$. Let $X$ and $Y$ be the
subsets of $V(H)$ that correspond to $I$ and $J$ in $W_H$. Then
\[
\|W_H-W_{G'}\|_\square \ge \int_{I\times
J}(W_{G'}-W_{H})=\frac{e_{G'}(S,T)}{n^2}-\frac{e_H(X,Y)}{n^2m^2} =
\frac{1}{4}-\frac{e_H(X,Y)}{n^2m^2}.
\]
Now because of the random labeling of $V(H)$, the expectation of the
second term is $1/8$, and the value of the second term will be
arbitrarily highly concentrated around this value if $m$ is
sufficiently large. There are only $n!$ possible relabelings of $G$,
so with large probability the second term will be less than $3/20$
for all of these, which proves our claim.

\subsection{Proof of \eqref{SQUARE-SS}, \eqref{SQUARE-ST-DISJ} and
\eqref{SQUARE-ST-COMPL}}
The upper bound in each of the three equations is trivial, so we only need
to prove
the lower bounds.  To prove these, we introduce the notation
\[
W(S,T)=\int_{S\times T}
W(x_1,x_2)\,dx_1dx_2
\]
for the value of the ``cut between $S$ and $T$''.

{\bf Proof of the lower bound in \eqref{SQUARE-SS}}:
This follows easily by
inclusion-exclusion applied to $W(S\cup T,S\cup T)$.

{\bf Proof of the lower bound in \eqref{SQUARE-ST-DISJ}}: To this
end, we first approximate $W$ by a weighted finite graph without
loops.  Next, we discard each point in $S\setminus T$ and $T\setminus
S$ with probability $1/2$, and then we add each point in $S\cap T$
uniformly at random either to the remaining points in $S\setminus T$,
or to the remaining points in $T\setminus S$.  In expectation, the
weighted number of edges  between the resulting, disjoint sets is off
by a factor of $1/4$, giving the first bound in
\eqref{SQUARE-ST-DISJ}.

{\bf Proof of the lower bound in \eqref{SQUARE-ST-COMPL}}:
To prove this bound, we consider two disjoint sets $S$ and $T$
and their complement $R=[0,1]\setminus(S\cup T)$.  We then express
each of the three cuts $W(S,T\cup R)$, $W(T,S\cup R)$, and $W(S\cup
T,R)$, in terms of $W(S,T)$, $W(S,R)$ and $W(T,R)$.  Combining these
three relations leads to the desired bound in
\eqref{SQUARE-ST-COMPL}.

\end{document}